\documentstyle[12pt,amssymb]{amsart}

\newcommand{\half}{{1\over 2}}
\newcommand{\wt}{\widetilde}

\newcommand{\PP}{{\mathbb P}}

\newcommand{\C}{{\mathbb C}}
\newcommand{\Q}{{\mathbb Q}}
\newcommand{\Z}{{\mathbb Z}}
\renewcommand{\phi}{\varphi}

\newcommand{\bX}{{\overline {X}}}
\newcommand{\bV}{{\overline {V}}}
\newcommand{\bS}{{\overline {S}}}
\newcommand{\G}{{\Gamma}}
\newcommand{\D}{{\Delta}}
\newcommand{\NS}{{\operatorname{NS}}}

\newcommand{\sym}{{\operatorname{sym}}}
\newtheorem{thm}{Theorem}[section]
\newtheorem{cor}[thm]{Corollary}
\newtheorem{lem}[thm]{Lemma}
\newtheorem{prop}[thm]{Proposition}

\newcommand{\proof} {{\it Proof.} }

\newenvironment{ex}{\smallskip\noindent\refstepcounter{thm}{\bf Example
\arabic{section}.\arabic{thm}.}}{\smallskip}
\newenvironment{rem}{\smallskip\noindent\refstepcounter{thm}{\bf Remark
\arabic{section}.\arabic{thm}.}}{\smallskip}

\newenvironment{defin}{\smallskip\noindent{\it Definition:\/}}{\smallskip}

\title {Two classes of
hyperbolic surfaces in $\PP^3$}

\author{Bernard Shiffman}
\address{Department of Mathematics,
Johns Hopkins University, Baltimore, MD
21218, USA} 
\email{shiffman@@math.jhu.edu}

\author{Mikhail Zaidenberg} 
\address{Universit{\'e}
Grenoble I, Institut Fourier, UMR 5582 CNRS-UJF, BP 74,
38402 St.\ Martin
d'H{\`e}res c{\'e}dex, France} 
\email{zaidenbe@@ujf-grenoble.fr}

\thanks{Research of the first author
partially supported by  NSF grant \#DMS-9800479.\\
\mbox{\hspace{11pt}}{\it 1991 Mathematics Subject Classification}: 
14J70, 14H35, 14E25, 14C20, 14M99, 32H20.\\
\mbox{\hspace{11pt}}{\it Key words}: projective surface,
generic projection, Kobayashi hyperbolic, Carath\'eodory
hyperbolic, Liouville variety, algebraic curve, symmetric square} 

\date{November 24, 1998}

\begin{document}

\begin{abstract} We construct two classes of singular
Kobayashi hyperbolic surfaces in $\PP^3$. The first consists of generic
projections of the cartesian square $V = C \times C$ of a generic genus $g
\ge 2$ curve $C$ smoothly embedded in $\PP^5$. These surfaces have
C-hyperbolic normalizations;
we give some lower bounds for their
degrees and provide an example of degree 32.  The second class of examples of
hyperbolic surfaces in  $\PP^3$ is provided by  generic projections of the
symmetric square $V' = C_2$ of a generic genus $g \ge 3$ curve $C$. The
minimal degree of these surfaces is 16, but this time the normalizations are
not C-hyperbolic. \end{abstract}

\maketitle

\section*{Introduction}
In this paper we construct new examples of  
Kobayashi hyperbolic   
surfaces in $\PP^3$. Concrete examples of Kobayashi hyperbolic  smooth
surfaces in
$\PP^3$ of degrees starting with $11$ were given in \cite{BrGre, Na, MaNo,
ElGu, SiYe}. Moreover, recently Demailly and El Goul \cite{DeElGu},  and
also McQuillan
\cite{MQ}, showed that a very generic surface in
$\PP^3$ of degree $d\ge 42$, resp.\ $d \ge 36$, is  Kobayashi hyperbolic, 
making the first step towards solving  the Kobayashi problem on
hyperbolicity  of a generic hypersurface in $\PP^m$  
of high enough degree. 
Note that examples of smooth hyperbolic
hypersurfaces in $\PP^m$ (for any $m \ge 3$) were constructed 
in Masuda-Noguchi \cite{MaNo}, and that the set of all hyperbolic
hypersurfaces  is open in the Hausdorff topology on the space of
hypersurfaces \cite{Za 1}. 

Our first examples of hyperbolic surfaces in $\PP^3$ have C-hyperbolic
normalizations. A complex space $X$ is called {\it C-hyperbolic} 
if there exists a non-ramified 
holomorphic covering $Y \to X$ 
with $Y$ being a Carath\'eodory 
hyperbolic complex space, i.e. 
any two points of $Y$ can be separated by 
a bounded holomorphic function \cite{Kob}. 
Since a Carath\'eodory hyperbolic complex space $Y$ 
is also Kobayashi hyperbolic, the base $X$ 
is Kobayashi hyperbolic, too. 
Thus, C-hyperbolicity implies Kobayashi hyperbolicity. 

A degree $d \ge 4$ smooth curve in $\PP^2$ is C-hyperbolic. 
However, for $m \ge 3$, hypersurfaces in $\PP^m$  
are not C-hyperbolic. 
Indeed, by the Lefschetz
hyperplane section theorem, they are simply connected, 
and so, do not admit non-trivial coverings. 
Instead, we construct
singular surfaces in $\PP^3$
with C-hyperbolic normalization by starting with a projective
embedding 
$X \hookrightarrow \PP^N$ of a C-hyperbolic surface $X$ and then considering
its general linear projection $\bX$ to $\PP^3$. 
It is classically known (see e.g.,
\cite{Pi}) that the normalization of $\bX$ is smooth and so,
coincides with $X$. 
Thus, the singular surfaces of type $\bX$ provide
examples of surfaces in $\PP^3$ with C-hyperbolic normalization.
In particular,
$X$ can be  a product $\G_1 \times \G_2$ of smooth projective curves
$\G_1,\,\G_2$ of genera $g_1,g_2\ge 2$; having the unit bidisc  
as its universal covering, the variety 
$X$ is C-hyperbolic.   However, a singular projective
surface $\bX$ with hyperbolic normalization $X$
can be non-hyperbolic even if $X$ has the unit bidisc as universal
cover; see \cite{KaZa} for 
an example. 

In \S \ref{sec-square},
we prove that for the cartesian square $V  = C \times C$
of a generic genus $g \ge 2$ curve $C,$ the generic projection
$\bV\subset\PP^3$ is, indeed,  Kobayashi hyperbolic (Theorem
\ref{hyp}). To justify our interest in the singular surface $ \bV $,  recall
that a hyperbolic hypersurface in $\PP^m$ possesses a hyperbolic
neighborhood, and thus all  hypersurfaces sufficiently close to $ \bV $ are 
hyperbolic, too  \cite{Za 1}.
In \S \ref{sec-symmetric}, we shift our attention to the symmetric
square $V'=C_2$ of a curve $C$ of genus $\ge 3$, and we show that
$V'$ is Kobayashi hyperbolic if $C$ is neither hyperelliptic nor
bielliptic (Proposition
\ref{hypersym3}). Using the symmetric square $V'$ rather than the
cartesian square $V$ as above, we produce an example of a Kobayashi
hyperbolic surface
$\bV'\subset \PP^3$ of degree $16$ (Theorem \ref{deg16}), which is
the smallest possible degree of $\bV'$ when $C$ has general
moduli\footnote{As usual, we say that a genus $g$ curve $C$ has {\it 
general moduli}, or is {\it generic}, if in the moduli space ${\cal M}_g$ it
belongs to the  complement of a certain countable union of proper
subvarieties.} (see below).  We observe that this time the normalization of
our surface fails to be C-hyperbolic; indeed, the universal cover of $C_2$
carries no nonconstant bounded holomorphic functions (Corollary
\ref{super-L}).

We would like to determine the minimal possible
degree of the surfaces $\bV$, resp.\ $\bV'$. With this aim in
mind, we let $\delta(Y)$ denote the minimal degree of the projective
embeddings of a projective variety $Y$.  In terms of this
notation, we want to know the smallest value of $\delta(C\times C)$, resp.\
$\delta(C_2)$.  In \S \ref{sec-square}, resp.\ \S\ref{sec-symmetric}, we
obtain some lower bounds for these minimal degrees  in
terms of the genus $g$ of $C$ (Proposition \ref{est}, resp.\
Proposition \ref{estsym}).
In particular, we show that (under the assumption of
generality of $C$)  $18 \le \delta(C\times C) \le 46$ if $g=2$, $20
\le \delta(C\times C) \le 32$ and $\delta(C_2)=16$ if $g=3$, and $20\le
\delta(C_2)\le 36$ if $g=4$ (see
Proposition \ref{est}, Example \ref{genus3}, Corollary \ref{min}(a),
Proposition \ref{estsym}, and Corollary \ref{newcor}).

In \S \ref{sec-divisors}, we investigate the divisors on $C\times C$ and
on $C_2$.  In
particular, we describe some classes of very ample divisors on symmetric and
cartesian squares of curves of genera 2, 3, and 4 (see
Theorems \ref{g2ample}, \ref{veryample3} and \ref{g2ample'}). The
generic projections of the projective embeddings given by these divisors
provide specific examples of the hyperbolic surfaces $\bV$ and
$\bV'$ described in \S\S \ref{sec-square}--\ref{sec-symmetric}.

In the course of preparation of this paper
we had  interesting and useful discussions 
with E. Amerik, F. Catanese, C. Ciliberto, B. Hassett, C. Peters and
V. Shokurov, and we would like to thank all of them.  
The second author also thanks the Johns Hopkins University at
Baltimore  and the Max-Planck-Institut f\"ur Mathematik 
at Bonn for their hospitality.
 
\bigskip
\section{Generic projection of a smooth surface into
$\PP^3$}\label{sec-projection}

Let $V$ be a smooth projective surface in $\PP^N$, $ N\ge 4$,
and let $\bV$ be its image
under a generic linear projection to $\PP^3.$ 
By the Severi Theorem
\cite{Sev} (see also \cite[(4.6)]{Mo, Pi, GriHarr}), the surface $\bV$ 
has only ordinary singularities. That is, 
$\bS := {\rm sing}\,\bV$ 
is a double curve of $\bV,$ i.e. a generic point 
$P \in \bS$ on $\bV$ is a point of 
transversal intersection of two smooth surface germs. 
The singularities of the curve $\bS$ itself 
are ordinary triple points, 
$P_1,\dots,P_t,$ say, which are also ordinary 
triple points of the surface $\bV,$ i.e. points of 
transversal intersection of three smooth surface germs.
Besides, there is also a certain number $Q_1,\dots,Q_p \in \bS$ 
of pinch points of the surface $\bV.$ They 
are smooth points of the double curve $\bS,$ 
and the local equation of 
$\bV$ around the point $Q_i$ is that of the Whitney umbrella
$x^2y - z^2 = 0.$ 

We will assume that  the
surface $V \subset \PP^N$ is not contained in a hyperplane and 
does  not coincide, up to projective transformations, 
with the image of $\PP^2$ under
the second Veronese 
embedding $\PP^2 \hookrightarrow \PP^5.$ 
By results of Severi, Moishezon and Mumford 
(see \cite[pp.~60, 72, 114-115]{Sev, Mo}), 
the double curve $\bS$ is irreducible, 
as well as its preimage
$S$ in the surface $V.$
The curve $S$ has $3t$ ordinary double points 
over the triple points $P_1,\dots,P_t$ of $\bS$ 
as the only singular points. 
The restriction 
$\pi\,|\,S\,:\,S \to \bS$ of the projection 
$\pi\,:\,V \to \bV$ 
(which is, at the same time, the normalization of $\bV$)
is generically $2 : 1$; it is ramified only over the pinch points 
$Q_1,\dots,Q_p \in \bS,$
and is an immersion at each of the double points of the curve $S$  (i.e.,
$\pi\,|\,S$ maps a neighborhood of each double point of $S$ injectively to a
neighborhood of a triple point of $\bS$).

The following simple observation reduces 
the Kobayashi hyperbolicity of the surface 
$\bV$ to those of the double curve $\bS.$

\medskip

\begin{prop} \label{khyp0} Let $\bX$ be a reduced compact complex space, 
and let 
$\pi\,:\,X \to \bX$ 
be the normalization of $\bX.$ Assume that the space $X$ is 
Kobayashi hyperbolic and let $S \subset X$, resp.\ 
$\bS:= \pi(S) \subset \bX$, be the
ramification  divisor, resp.\ the branching divisor, so that  
the restriction $\pi\,|\,(X \setminus S)\,:\,X \setminus S
\to \bX \setminus \bS$ is biholomorphic. 
Then $\bX$ is Kobayashi hyperbolic
if and only if  $\bS$ is Kobayashi hyperbolic.\end{prop}

\medskip

\noindent \proof Clearly, if the space $\bX$ is hyperbolic, 
then so is the subspace $\bS$ of $\bX.$
By the Brody Theorem \cite{Br} 
(see \cite{Ki}, \cite[(3.6.3)]{Kob} or \cite{Za 3} for the case of 
complex spaces), the compact complex space 
$\bX$ is hyperbolic iff any holomorphic mapping 
$f\,:\,\C \to \bX$ is constant. Assuming that 
the branching divisor $\bS$
is hyperbolic, we may restrict the consideration to 
the mappings $f\,:\,\C \to \bX$ with 
the image not contained in $\bS.$ In this case
$f$ can be lifted to $X$ (see \cite{PeRem}) and hence, 
in virtue 
of the hyperbolicity of the complex space
$X,$ it is constant. \qed

\medskip Applying Proposition~\ref{khyp0} to the case where
$V \subset \PP^N,\ N\ge 4$, is a smooth Kobayashi 
hyperbolic surface 
(e.g., $V$ can be isomorphic to
the cartesian product 
of two smooth projective curves $\G_1$ and $\G_2$ of genera 
$g_1,\,g_2 \ge 2$) and
$\bV \subset \PP^3$ is a generic projection of $V$
with the irreducible double curve $\bS,$ we obtain the following
statement:

\medskip

\begin{cor} \label{khyp} {Let $\bV \subset \PP^3$ be a generic projection of
a Kobayashi 
hyperbolic  smooth projective surface.  Then $\bV$ is Kobayashi hyperbolic 
iff  the double curve $\bS$ is Kobayashi hyperbolic; i.e., 
iff the geometric genus of $\bS$ is at least $2.$} \end{cor}

\medskip

We denote by $H$ a generic hyperplane section 
of $V$ in $\PP^N$
regarded as a very ample divisor on $V.$
Let $n = H^2$ be 
the degree of $V,$ and $g_H$ be the genus 
of the smooth curve $H$ on $V,$ i.e. the sectional genus
of the surface $V$ in  $\PP^N.$
Let $c_1, \, c_2$ 
be the Chern classes of $V$ and 
$K = K_V$ be the canonical divisor of $V.$ 
Denote by $d_{\bS}$, resp.\ $g_{\bS},$ 
the degree, resp.\ 
the geometric genus,
of the double 
curve $\bS;$ $t$, resp.\ $p$, denotes, as above, 
the number of triple points, resp.\ 
pinch points, of $\bV.$ We also denote by $b$ 
the number of double points 
of a generic projection $V \to \PP^4.$

We have the  following relations \cite[(4.6)]{GriHarr} \cite{PetSi}  \cite{Pi}
\cite[IX.5.21]{SemRo} (cf. also \cite[(8.2.1)]{BeSo},
\cite{Mo} and \cite[A4]{Hart}): 
\begin{eqnarray} \label{c1} c_1^2 &=& n(n-4)^2 - (3n-16)d_{\bS}
+ 3t - b \\ \label{c2} c_2 &=& n(n^2-4n +6) - (3n-8)d_{\bS} + 3t - 2b \\
\label{chi} \chi({\cal O}_V)= {c_1^2 + c_2 \over 12} &=&
{n-1 \choose 3} - d_{\bS}(n-4) + g_{\bS} + 2t \\
\label{p1} p &=& 2d_{\bS}(n-4) - 6t - 4(g_{\bS}-1) \\
\label{b1} b &=&{1\over 2} [n(n - 10) - 5H\cdot K - c_1^2 + c_2]\nonumber\\
&=&{1\over 2} [n(n - 5) - 10(g_H - 1)- c_1^2 + c_2]
\end{eqnarray}

(Note that equations (\ref{c1})--(\ref{p1}) satisfy a linear relation.)

\medskip

\begin{prop} In the notation as above, we have
\begin{equation}\label{p} p = c_1^2 - c_2 + 2n(n-5) - 8d_{\bS} =
c_1^2 - c_2 + 2n + 8(g_H - 1) 
= c_1^2 - c_2 + 6n + 4H\cdot K \,\end{equation}
where
\begin{equation}\label{d} d_{\bS}= {n-1 \choose 2} - g_H \,.\end{equation}
\end{prop}

\noindent \proof Subtracting (\ref{c2}) from (\ref{c1}) gives the first
equality  in (\ref{p}), the second one is provided by plugging in  the value
of $d$ from (\ref{d}),  and the third one by the Adjunction Formula

\begin{equation} \label{adj} 2g_H-2=H\cdot (H+K)=n+H\cdot K\,. \end{equation}
 
To prove (\ref{d}),
consider a generic hyperplane $h \subset \PP^N$ 
which contains the center of the projection
$\pi\,:\,\PP^N \to \PP^3.$ Then the image
$l:=\pi(h)$ is a generic plane in $\PP^3.$ 
Respectively,
the image $L:=\pi(H) = l\cdot \bV$ 
of the hyperplane section 
$H = h \cdot V$ is a generic plane section of the
surface $\bV = \pi(V) \subset \PP^3.$ Hence,
$L$ is a degree $n$ plane curve with $d_\bS$ nodes, 
which are smooth points of the double curve $\bS$, and no other
singularities. Now (\ref{d}) is the genus formula for the plane nodal
curve $L.$ 
\qed

\medskip

\begin{cor} {We have the following expressions for $g_{\bS}$ and $t:$}
\begin{eqnarray} \label{g} g_{\bS} &=& {1 \over 2}(n^2-7n +26) +
(n-12)g_H- {5c_1^2 - 3c_2 \over 4} \\ \label{t1} t &=&{1 \over
6}(n^2-3n+8)(n-6)  - (n-8)g_H + {2c_1^2 - c_2 \over 3} \end{eqnarray}\end{cor}

\noindent \proof
{From} (\ref{p1}) and (\ref{p}), resp.\ (\ref{chi}),
we obtain the equations
\begin{eqnarray} \label{s1} 3t + 2g_{\bS} &=& (n-4)d_{\bS} - p/2 + 2 = 
nd_{\bS}-n(n-5)-{c_1^2 - c_2 \over 2} \;, \\
\label{s2} 2t + g_{\bS} &=& (n-4)d_{\bS} 
- {n-1 \choose 3} + {c_1^2 + c_2 \over 12}\;, \end{eqnarray}
which together with (\ref{d}) give us (\ref{g})--(\ref{t1}). \qed

\bigskip
\section{Generic projection of a surface $V = C \times C$
into $\PP^3$}\label{sec-square}

{From} now on we let $V = C \times C$ where 
$C$ is a  
smooth projective curve of genus $g \ge 2.$  Fix a point $p_0\in C$ and
consider the horizontal, resp.\ vertical, divisor $E=C \times
\{p_0\}$, resp.\  $F=\{p_0\}\times C.$ The canonical divisor $K = K_V$ 
is numerically equivalent to the divisor 
$(2g-2)E + (2g-2)F$ on $V,$ and we have 
\begin{equation}\label{chern} c_1^2=K^2 = 2(2g-2)^2,\quad
c_2 = e(V) = e(C)^2 = (2g-2)^2 = c_1^2/2 \end{equation}
where $e$ stands for the Euler characteristic.
Hence, the signature $\tau = \tau(V)$ vanishes: 
$\tau = {1\over 3}(c_1^2 - 2c_2) = 0.$ 
By the Noether Formula for the holomorphic 
Euler characteristic,
we have
$$\chi({\cal O}_V) = 
{c_1^2 + c_2 \over 12} = (g-1)^2\,.$$
Furthermore, for the geometric genus resp.\ 
the irregularity of the surface $V$ we have 
$$p_g(V) = h^{0,2}(V) = g^2,\,\,\,\,\,\,\,\,\,
q(V) = h^{0, 1}(V)= 2g\,.$$
Denote by $\Delta$ the diagonal of 
the cartesian square $V = C \times C.$
Let $\Sigma$ be the subgroup of the 
Neron-Severi group $\NS(V)$
generated by the classes of 
$E,\,F,\,\Delta$ modulo 
numerical equivalence. By a theorem of 
Hurwitz \cite{Hu} (see also \cite[(2.5)]{GriHarr}), for a 
genus $g \ge 2$ curve $C$ with generic 
moduli we have 
$\Sigma = \NS(V)$. Set 
$D_{a,a',k} = (a+k)E + (a'+k)F - k\Delta.$
Using the fact that $\Delta^2 = 2 - 2g$ 
we obtain the following
standard formulas from the theory of correspondences
\cite{Fu, GriHarr, SemRo, We}: 
$$ D_{a,a',k}\cdot
D_{c,c',l} = ac' + a'c - 2gkl, \,\,\,D_{a,a',k}^2=2(aa'-gk^2)\,,$$
$$D_{a,a',k}\cdot E = a',\,\,\,D_{a,a',k}\cdot F = a\,.$$
(The pair $(a,\,a')$ is called the {\it bidegree} 
of $D:=D_{a,a',k},$ and $k$
is called the {\it valence} of $D.$)
Furthermore, for $a=a'$ we denote 
$D_{a,a',k}$ as $D_{a,k};$ thus,
$\Delta \equiv D_{1,-1}$ and $K \equiv D_{2g-2,0}$ where $\equiv$ 
stands for numerical equivalence. 
We have:
$$\Delta^2 = 2 - 2g,\,\,\,K^2 = 2(2g-2)^2\,,$$
$$D_{a,k}\cdot D_{a',l} = 2(aa' - gkl);
\,\,\,D_{a,k}^2 = 2(a^2 - gk^2)\,,$$
$$D_{a,a',k}\cdot \Delta = a+a'+2gk,
\,\,\,D_{a,a',k}\cdot K = 2(g-1)(a+a')\,,$$
in particular, 
$$D_{a,k}\cdot  \Delta = 2(a+gk),
\,\,\,D_{a,k}\cdot K =4(g-1)a\,.$$

\begin{prop} Let the notation be as in Section 1 above. 
For a very ample divisor $H \equiv D_{a,a',k} \in \NS(V)$ on the surface 
$V= C \times C$ we have:
\begin{eqnarray} \label{n} n &=&  2(aa'-gk^2) \\
\label{gH} g_H &=& {n \over 2}+1 + (g-1)(a+a') \\
\label{dS} d_{\bS} &=&  {1 \over 2}n(n-4)-(g-1)(a+a') \\
\label{gS} g_{\bS} &=& {1 \over 2}(2n^2 -17n +2) + 
(n-12)(g-1)(a+a')-7(g-1)^2 \\
\label{p2} p &=& 2\left[3n+4(g-1)(a+a') + 2(g-1)^2\right] \\
\label{tC} t &=& {1\over 6}n(n^2-12n+44)-
(n-8)(g-1)(a+a') +4(g-1)^2\\
\label{b2} b &=& {1\over 2}n(n-10) - 5(g-1)(a+a') - 2(g-1)^2  \nonumber\\
&=&{1\over 2}n(n-5) - 5(g_H-1)- 2(g-1)^2 
\end{eqnarray} \end{prop}

\noindent \proof Equation (\ref{n}) follows from the fact that
$n=\deg\,V =H^2.$   By the
Adjunction Formula (\ref{adj}), 
$$2g_H-2= H\cdot K +n = 2(g-1)(a+a')+n\,,$$ 
and (\ref{gH}) follows.  
Substituting (\ref{chern}) and (\ref{gH}) into
(\ref{b1})--(\ref{d}), we obtain (\ref{dS})--(\ref{b2}). \qed

\medskip

We use these formulas to prove the following inequality:

\medskip

\begin{prop} \label{est} Let $C$ be a genus $g\ge 2$ curve with generic
moduli. The minimal degree $n=\delta(C\times C)$ 
of a projective embedding 
$C\times C\hookrightarrow \PP^N$
satisfies the inequality \begin{equation}\label{in} n(n-10) \ge
4(g-1)(g-1+5\delta(C))\,\end{equation}  In particular, $\delta(C\times C)
\geq 18$ for $g=2$, $\delta(C\times C) \geq 20$ for $g=3$, $\delta(C\times
C) \ge 2g+16$ for $g\ge 4$, and $\delta(C\times C)\ge 2(g-1) +5\sqrt{2g}$
for  $g \gg 0.$ \end{prop}

\noindent \proof Let $H\equiv D_{a,a',k}$ be a hyperplane section of
the minimal degree embedding $V=C\times C\hookrightarrow \PP^N$. Since 
$a=H\cdot F =\deg (F\hookrightarrow
\PP^N) \ge \delta(C)$ and similarly $a'=H\cdot E\ge \delta(C)$, it follows from
(\ref{gH}) that \begin{equation} \label{gHin} g_H\ge \half(n+2)
+2\delta(C)(g-1)\,.\end{equation} On the other hand, from
(\ref{b2}) we have the inequality
$$2b=n(n-5) - 10(g_H-1) -4(g-1)^2 \ge 0\,,$$
or equivalently
\begin{equation}\label{in2} g_H \le {n^2 - 5n + 10 - 4(g-1)^2\over 10}\,.
\end{equation}
Combining (\ref{gHin}) and (\ref{in2}), we obtain (\ref{in}). 

We recall the standard lower bounds for $\delta(C)$.  If $g=2$, then
$\delta(C)=5$ since the genus of a smooth curve  in $\PP^3$ of degree
$\le 4$ is at most 1 unless it is a plane curve of genus 3 (and degree 4). 
In general, the geometric genus $g$
of a projective curve of degree $\delta$ satisfies the inequality 
$g\le \half
(\delta-1)(\delta-2)$; thus
\begin{equation}\label{delta}\delta(C)\ge {3\over
2}+\sqrt{2g+{\textstyle{1\over 4}}}\ \ {\rm for }\ \ g \ge 3\;,\quad 
\delta(C)= 5\ \ {\rm for }\ \ g=2\;.\end{equation}  

The bounds for $g=2,\,g=3,\,g\ge 4$ follow from (\ref{in}), the fact that
$n$ is even, and the estimate $\delta(C)\ge 5$ except for $g=3$ in which
case $\delta(C)\ge 4$.  To obtain the estimate for large $g$, we
note that by (\ref{in}) and (\ref{delta}), \begin{eqnarray*}n(n-10)
&\ge& 4(g-1)\left(g  +5\sqrt{2g+{\textstyle{1\over 4}}}+{13\over 2} \right)\\
&=&4g^2\left(1-g^{-1}\right)\left(1+5\sqrt{2}\,g^{-\half} +{13\over 2} g^{-1}
+ O(g^{-{3\over 2}}) \right)\,. \end{eqnarray*} Therefore $$\left({n-5\over
2g}\right)^2 \ge 1+ 5\sqrt{2}\, g^{-\half} + {11\over 2} g^{-1}
+O(g^{-{3\over 2}})\,,$$ which gives us  $$n-5 \ge 2g +5\sqrt{2}\,g^\half -7
+O(g^{-\half})\,.$$ Since $n$ is an integer, the large $g$ estimate follows.
\qed

\medskip

\begin{rem} \label{rm2.3}
For generic $C$ of genus $g \ge 4,$ we have values of
$\delta(C)$ larger than those given by (\ref{delta}) (see \cite{Hart},
\cite{GriHarr}).  For instance, $\delta(C) = 6$ for a non-hyperelliptic genus
$4$ curve $C$ \cite[p.40, (D-4, D-7)]{ACGH}.  These higher bounds in turn
give higher lower bounds for $n.$ Anyhow, the bounds for $\delta(C\times C)$
obtained from (\ref{in}) are probably far from being sharp.

Being a constructive integer-valued function on the moduli space  ${\cal
M}_g$ of genus $g$ curves, in general, $\delta(C)$ is not semi-continuous,
i.e. under specialization it can  increase as well as decrease. For example,
if $g=3$ then $\delta(C)=4$ outside of the locus of hyperelliptic curves,
whereas by Halphen's Theorem \cite[IV.6.1]{Hart}, $\delta(C)=6$ if $C$ is
hyperelliptic.  On the other hand, if $g=6$ then $\delta(C)\ge 8$ at a
generic point  $C \in {\cal M}_6,$ $\delta(C) = 7$ precisely at the locus
$L_1$ of smooth curves of  type $(3,\,4)$ on a smooth quadric $Q \simeq
\PP^1 \times \PP^1$ in  $\PP^3,$ and $\delta(C) = 5$ precisely at the locus
$L_2$ of smooth plane quintics (see \cite[IV.6.4.2]{Hart}); we have $\dim
{\cal M}_6 = 15,\,\,\,\dim L_1 = 11$ and $\dim L_2 = 13.$ \end{rem}

\medskip 

\begin{cor} \label{cin} \ Let $V=C\times C\hookrightarrow \PP^N$ be as in
Proposition~\ref{est}.  Then $g_\bS \ge 225$ for $g=2$, $g_\bS\ge 331$ for
$g=3$, and $$g_\bS \ge 17g^2+81g+74\quad \mbox{for}\quad g\ge 4\,.$$ In fact,
\begin{equation} \label{cinlargeg} g_\bS > 4\sqrt{2}g^{{5\over 2}}\quad
\mbox{for}\quad g \gg 0\,.\end{equation}
\end{cor}

\noindent \proof It follows from
(\ref{gS}) that 
\begin{equation}\label{in3} g_{\bS} \geq {1 \over 2}(2n^2
-17n +2) + 2\delta(C)(n-12)(g-1)-7(g-1)^2 \;.\end{equation} The result for
$g=2,\,g=3,\,g\ge 4$, resp., follows by substituting the bounds $n\ge 18,\
20,\ 2g+16$, resp., and $\delta(C)\ge 5,\ 4,\ 5$, resp., into
(\ref{in3}).  To verify (\ref{cinlargeg}), we substitute $n\ge 2g
+5\sqrt{2g}-2$ and (\ref{delta}) into (\ref{in3}) to obtain
$$g_\bS\ge 4\sqrt{2}g^{{5\over 2}} +17g^2+O(g^{{3\over 2}})\,.$$
\qed

\medskip

Corollaries \ref{khyp} and \ref{cin} lead to the following conclusion.

\medskip

\begin{thm} \label{hyp} {Let $V = C \times C$ where $C$ 
is a smooth projective curve of genus $g \ge 2.$ 
Then a generic projection $\bV \subset \PP^3$ of 
the image of any projective embedding $V \hookrightarrow \PP^N$ 
is a Kobayashi hyperbolic surface with C-hyperbolic
normalization.}\end{thm} \smallskip

\begin{ex} \label{genus3} 
Let $C$ be a non-hyperelliptic smooth projective
curve of genus $g = 3.$ It is well known that 
the canonical divisor $K_C$ of degree $2g-2 = 4$ 
is very ample,
and it provides the canonical embedding 
$C \hookrightarrow \PP^2$
onto a smooth plane quartic. In turn, 
the canonical divisor 
$K = K_V \equiv 4E + 4F \equiv D_{4,0}$ 
of the cartesian square $V = C \times C$ yields the 
Segre embedding $V \hookrightarrow \PP^8$
onto a smooth surface of degree $n = K^2 = 32.$ \end{ex}

In the case where $C$ has genus 2, any effective divisor of degree 5
gives an embedding $C \hookrightarrow \PP^3.$ Thus, 
the divisor $D_{5,0}$ on
the cartesian square $V = C \times C$ yields the embedding 
$V \hookrightarrow
\PP^{15}$ onto a smooth surface of degree $n = 50.$
 
In each of these two cases, formulas (\ref{n})--(\ref{b2}) 
give us the following numerical data for a generic projection 
of $V=C\times C$ to $\PP^3$:
$$\begin{array}{l@{\,}l@{\,}l@{\,}l@{\,}l@{\,}l@{\,}l@{\,}l@{\,}l@{\,}l}
g=2,&H\equiv
D_{5,0}&\Rightarrow &n=50,&g_H = 36,&d_\bS=1140,&g_{\bS} = 2449,&p = 
384,&t=15784,&b=948\\ g=3,&H\equiv D_{4,0} &\Rightarrow &n=32,&g_H =
33,&d_\bS=432,&g_{\bS} = 1045, &p =  336,&t=3280,&b=264 \end{array}$$

\begin{rem} \label{rm2.1} A very ample divisor $H \equiv D_{a,a',k} \in
\NS(V)$ defines an embedding $V \hookrightarrow \PP^4$ iff $b=0,$ i.e. in
view of (\ref{b2}), iff \begin{equation}\label{P4} n(n-10) = 10(g-1)(a+a') +
4(g-1)^2\,. \end{equation} But we don't know whether there is a smooth
surface in $\PP^4$ isomorphic to $V=C\times C$.\footnote{For $g=2,$ the
first integer solution of (\ref{P4}) which also satisfies the necessary
conditions for very ampleness given in Remark \ref{(e)} below is
$a=511,\ a'=79,\ |k|=142$; but the divisor $D_{511,79,\pm 142}$ is probably
not very ample.} Indeed Van de Ven conjectured that the irregularity of a
smooth surface in $\PP^4$ can be at most 2 \cite[Problem 8]{OkVdV}, which
would imply that $V$ cannot be embedded in $\PP^4$ and thus $b>0$. Any
better estimate of $b$ from below would lead to a better lower estimate for
$\delta(C\times C)$.

The next proposition shows the nonexistence of such
an embedding at least in certain cases (cf.\ \cite{CGB, IdMe}).  \end{rem}

\medskip

\begin{prop} \label{emb}
Let $C$ be a generic curve of genus $g \ge 2,$ and
let $H \equiv D_{a,a',k}$ be a very ample divisor on the cartesian square
$V=C\times C.$ Then $H$ cannot provide an embedding 
$V \hookrightarrow \PP^4$ in each of the following cases:

\smallskip\noindent {\rm (a)} $H$ is a non-special divisor
(i.e., $h^1({\cal O}_V(H))=0$) and $g \le 13;$

\smallskip\noindent{\rm (b)}  $H$ is a non-special divisor and $a+a' \ge
g+1$;

\smallskip\noindent{\rm (c)}  $a=a'$ and $g=2.$ \end{prop}

\noindent \proof
(a) Suppose $H$ were to give an embedding $V \hookrightarrow \PP^4.$ By
Severi's Theorem \cite{Sev} (see also \cite[ex.~IV.3.11.b]{Hart}), 
a smooth surface $X$ in $\PP^4$ is linearly
normal, i.e. $h^0({\cal O}_X(1))=5,$  unless  it is a projection of the second
Veronese embedding  $\PP^2  \hookrightarrow \PP^5$.  Hence $h^0({\cal
O}_V(H))=5$, and thus by the Riemann-Roch and Adjunction Formulas together
with (\ref{gH}) and (\ref{P4}), we obtain: \begin{eqnarray} \chi({\cal O}(H))
&=&  5- h^1({\cal O}(H)) +h^2({\cal O}(H)) ={\half H\cdot (H-K)} + {c_1^2 +
c_2 \over 12}\nonumber\\&=&n-g_H+1+(g-1)^2 = \half n-(g-1)(a+a') 
+(g-1)^2\nonumber\\ &=&{15n-n^2+ 14(g-1)^2 \over 10} \,. \label{2P4}
\end{eqnarray}
 
Assuming that the embedding $V \hookrightarrow \PP^4$
is non-special, i.e. $h^1({\cal O}(H))=0$ (cf. \cite{IdMe}), 
from (\ref{2P4}) we get the inequality 
\begin{equation}\label{3P4} (n-5)(n-10)\leq 14(g-1)^2 \,. \end{equation}
On the other hand, from (\ref{in})
and (\ref{delta}), we obtain:

\bigskip
\hfil\begin{tabular}{|r|c|c|c|r|c|c|}\cline{1-3} \cline{5-7}
$g$ & $\delta(C)\ge$ & $n\ge$ & \quad\quad\quad\quad &$g$ & $\delta(C)\ge$ &
$n\ge$\\ \cline{1-3}\cline{5-7} 2&5&18&&8&6&38\\ 
\cline{1-3}\cline{5-7}3&4&20&&9&6&40\\ 
\cline{1-3}\cline{5-7}4&5&24&&10&6&44\\
\cline{1-3}\cline{5-7}5&5&28&&11&7&48\\
\cline{1-3}\cline{5-7}6&5&30&&12&7&50\\
\cline{1-3}\cline{5-7}7&6&36&&13&7&54\\ \cline{1-3}\cline{5-7}
\end{tabular}\hfil\bigskip 

The lower bounds for $n$ from this table contradict (\ref{3P4}). 
This proves (a).

\smallskip

(b) If $H$ is an embedding and (b) were to hold, then by (\ref{P4}) and
(\ref{3P4}), we would have
\begin{equation} \label{est'}
5(n-10)=n(n-10)-(n-5)(n-10)\ge 10(g-1)(a+a'-g+1)\ge 20(g-1)\end{equation}
and thus by (\ref{3P4}) and (\ref{est'}),
$$14(g-1)^2\ge (n-10)^2 \ge 16(g-1)^2\;,$$ 
a contradiction. 

\smallskip

For the proof of (c) see Corollary \ref{min}(d) below.
\qed

\bigskip
\section{The symmetric square of a curve of genus $\ge
3$}\label{sec-symmetric}

As in \S \ref{sec-square}, we let $C$ be a curve of genus $g$,
except we now assume that $g\ge 3$.  We  still let
$V=C\times C$ denote the cartesian square, and we let 
$V':=C_2=C\times C/\{1,\sigma\}$
denote the symmetric square, where
$\sigma:V\to V$ is the involution $\sigma(z,w)=(w,z).$ 
More generally, consider the $d-$th symmetric power $C_d,\,\,\,d>0,$ 
of $C.$ Recall that $C_d$ can
be identified with the space of effective divisors of degree $d$ on $C,$ 
and the fibres of the Abel-Jacobi morphism $u_d\,:\,C_d \to J_C$ 
into the Jacobian variety $J_C$ 
can be identified with the complete degree $d$ linear systems on $C.$ 

By a {\it Brody curve} in a compact hermitian complex manifold $M$  we mean
a non-constant holomorphic map $f\,:\,\C \to M$ with bounded derivative
$f'\,:\,\C \to TM$ (see e.g. \cite{Za 3}).  Recall that by the
Brody Theorem \cite{Br}, $M$ is hyperbolic iff it does not admit any Brody
curve.  In the next proposition, we use the same approach as in
M.~Green \cite{Gre} to describe Brody
curves in the symmetric powers $C_d$ of a smooth projective curve 
$C$. 

\medskip \begin{prop}\label{brodycurve} Let $f\,:\,\C \to C_d$ be a Brody
curve. Then either

\medskip\noindent {\rm (a)} $u_d \circ f (\C) = $const, and then there exists a
linear pencil $g_d^1$ on $C,$ or

\medskip\noindent {\rm (b)}  the image $u_d \circ f(\C)$ lies on a smooth
abelian subvariety of $W_d := u_d(C_d) \subset J_C.$ \end{prop}

\noindent \proof (a) If $u_d \circ f(\C) = $const, i.e. $f(\C)$ is contained in
a  fiber of $u_d,$ then this fiber represents a complete linear system
$g_d^r$  on $C$ of positive dimension $r;$ in particular, it contains a
linear pencil $g_d^1.$

(b) Otherwise, $u_d \circ f\,:\,\C \to J_C$ is a Brody curve, and so is the
lift ${\widetilde f}\,:\,\C \to \C^g.$ As noted in \cite{Gre}, the
derivative ${\widetilde f}'$, being bounded, must be constant; that is,
${\widetilde f}$ is an embedding onto an affine line in $\C^g$, and the
closure   of the image $u_d \circ f(\C)$ in $J_C$ is a shifted
subtorus contained in the subvariety $W_d = u_d(C_d)$.  \qed

\medskip \begin{cor} For a generic curve $C$  of genus 
$g \ge 3$ and for any 
$d < g/2 + 1$ the symmetric power $C_d$ is Kobayashi hyperbolic,
while for any genus $g$  curve $C$ and any $d \ge g/2 + 1$ 
it is not hyperbolic. 
\end{cor}

\medskip

\noindent \proof By Theorems 1.1 and 1.5 in \cite[Ch. V]{ACGH}, for $d
\ge {g \over 2} + 1$ any genus $g$ curve $C$ possesses a linear pencil
$g_d^1$ (indeed, for the Brill-Noether number $\rho=\rho(g,\,d,\,1)=2d-g-2,$
we have $\rho \ge 0$). Hence, the symmetric power $C_d$ contains a smooth
rational curve, and so, it fails being Kobayashi hyperbolic. While for $d <
{g \over 2} + 1$ (i.e. for $\rho < 0$) a genus $g$ curve $C$ with generic
moduli has no $g_d^1.$ Furthermore, its Jacobian variety is simple, i.e. it
does not contain proper abelian subvarieties.  Indeed, it is known that in
the moduli space ${\cal M}_g$ of genus $g$ curves, the locus of curves with
non-simple Jacobian varieties is a countable union of subvarieties of
codimension at least $g-1$; see \cite[(3.4)-(3.7)]{CGB}. Hence by
Proposition \ref{brodycurve} and Brody's Theorem, the variety $C_d$, having
no Brody curve, must be Kobayashi hyperbolic. \qed

\medskip In the case of the symmetric square $V=C_2$ of a genus $g$ curve
$C,$  the next proposition provides more precise information.
Recall that
$V'\simeq \PP^2$ if $g=0,$ 
$V'$ is a ruled surface over $C$ if $g=1$, and $V'$ is the Jacobian $J_C$ 
blown up at a point if $g=2,$ so that the symmetric square of a genus
$g\le 2$ curve cannot be hyperbolic.

\medskip \begin{prop}\label{hypersym3} 
Let $C$ be a smooth projective
curve.  Then the following are equivalent:

\smallskip\noindent {\rm (i)} the surface $C_2$ 
is Kobayashi hyperbolic;
 
\smallskip\noindent {\rm (ii)} $C_2$ does not 
contain any rational or
elliptic curves;

\smallskip\noindent {\rm (iii)} the curve $C$ 
is neither hyperelliptic nor
bielliptic.
\end{prop}

\medskip Proposition \ref{hypersym3} is an immediate consequence of
Brody's Theorem, Proposition \ref{brodycurve} and the following two lemmas,
the first of which is well known and the second due to Abramovich and
Harris
\cite[Th.~3]{AH}.

\medskip

\begin{lem}\label{gamma}  The symmetric square  $V' = C_2$ contains a
rational curve $\G$ iff the curve 
$C$ is hyperelliptic.
\end{lem}

\medskip

\noindent \proof  
Since the Jacobian variety $J_C$ has no rational curve, 
$\G$ should coincide with a fibre of the Abel-Jacobi map 
$u_2\,:\,C_2 \to J_C,$ thus representing a pencil 
$g^1_2$ on $C$. Therefore, $C$ is a hyperelliptic curve.
The converse implication is easy.
\qed

\medskip

\begin{lem}{\bf (Abramovich-Harris)} \label{gamma-AH}  The
symmetric square  $V' = C_2$ of a genus 
$g\ge 3$ curve $C$ contains an
elliptic curve $\G$ iff the curve 
$C$ is bielliptic, 
i.e. there exists a $2:1$ morphism 
$f\,:\,C \to T$ of $C$ onto a smooth elliptic curve $T$.  
\end{lem}

\noindent\proof \cite[Th.~3]{AH}. \qed
\medskip

\begin{rem} \label{ell} (a) We easily see that the curve $\G$ in Lemma
\ref{gamma} is given as  $\G=f^{\vee} (\PP^1)$, where $f:C\to\PP^1$ is a 
$2:1$ morphism and $f^{\vee}\,:\,\PP^1 \ni z \longmapsto f^*(z) \in C_2$. It
can also be shown using
\cite[(3.2)]{AP} that for $g\ge 4$, the curve $\G$ in Lemma \ref{gamma-AH}
is similarly given as  $\G=f^{\vee} (T)$.
In Proposition \ref{gammasq} below, we show that for a curve 
$\G$ as above, $\G^2=1-g <0$. 

Observe that for a genus 2 curve $C$, 
the symmetric square 
$V'=C_2$ contains an elliptic curve iff $V'$ is an elliptic surface
(and hence, iff the Jacobian $J_C$ 
is isogenous to a product of two elliptic curves). Thus, if $g=2$ 
there may exist
smooth elliptic curves $\G$ on $C_2$ with $\G^2=0$ and so these curves
are not of the form $f^{\vee} (T)$.

\medskip

(b) Since the hyperelliptic involution is unique,  there can be only one
rational curve on $V'$, for $g\ge 2$. In contrast, if $C$ is bielliptic and
$2\le g\le 5$, there  can be several elliptic curves on $V'$. 
For instance, the Fermat quartic $C=\{x^4+y^4+z^4=0\}$ in 
$\PP^2$ admits 15 different involutions 
($(x:y:z) \longmapsto (-x:y:z)$, 
$(x:y:z) \longmapsto (\zeta y:\zeta^{-1}x:z),\,\,\,\zeta^4=1,$ etc.),
which are elements of order 2 in the automorphism group 
Aut$\,C\simeq (\Z/4\Z)^2 \rtimes S_3$ (see \cite[pp. 274--275]{KuKo}) and
which provide 15 different elliptic curves on $V'=C_2$.

Another example is the Klein quartic 
$C=\{xy^3+yz^3+x^3z=0\}$ in $\PP^2$ \cite[Ch. 8]{Kl}. 
The automorphism group Aut$\,C \simeq$ PSL$_3(\Z/2\Z)\simeq$ PSL$_2(\Z/7\Z)$ being a simple group of order 168 (which is the maximal possible one
for a genus 3 curve) \cite{Ku}, 
it is generated by 21 reflections, i.e. elements of order 2,
which form a conjugacy class. Each of them defines a bielliptic involution on $C$ and whence, an elliptic curve of the above type on the surface $V'=C_2$.
 
See also \cite[VI.F12-F14]{ACGH} 
for an example of a genus 5 curve bielliptic in ten different ways, 
so that $V'=C_2$ possesses 10 different elliptic curves.
However \cite[VIII.C-2]{ACGH}, a genus $g\ge 6$ curve $C$ can have only one
bielliptic structure,  and hence $V'=C_2$ may possess only one elliptic curve
$\G$.   

\medskip

(c) As was noted in \cite[Th. 3]{AH} (see also \cite[VIII.C-1]{ACGH}), a
bielliptic hyperelliptic curve $C$ is a type $(2,\,4)$ curve on a smooth
quadric,  and therefore of genus at most 3. (Indeed, a birational embedding 
$C \hookrightarrow \PP^1\times \PP^1$ is given by the pencils  $g^1_2$ and
$g^1_4$ on $C$ where the first one is the hyperelliptic one, and the second
one is provided by a hyperelliptic pencil $g^1_2$ on an elliptic curve
$T$ under a $2:1$ morphism $C\to  T$.) For instance
\cite[I.H-6]{ACGH}, $C=\{y^2z^4-x^6+z^6=0\}\subset \PP^2$ is such a curve of genus 2.

If $R$ is a genus 2 curve, then to any index 2 subgroup of the 
fundamental group $\pi_1(R)$ there corresponds a nonramified Galois $2:1$ covering $C \to R$
where $C$ is a genus 3 curve. 
Let $\sigma \in $Aut$\,C$ be the generator of the Galois group 
$\Z/2\Z$. Since $\sigma$ acts freely on $C$, and 
any involution of a smooth plane quartic has fixed points, 
$C$ is hyperelliptic; denote by $i$ the hyperelliptic involution on $C$. 
Then the involution $\sigma':=\sigma\circ i = i\circ\sigma \in $Aut$\,C$ 
has 4 fixed points (that is, the union of two common orbits of 
$i$ and $\sigma$; they come from the two fixed points 
of the induced action of $\sigma$ on the canonical model of $C$). 
Thus $\sigma'$ is a bielliptic involution. 
Therefore (Farkas-Accola; see \cite[Lemma 8]{KuKo}), 
any genus 3 curve which dominates 
a genus 2 curve is hyperelliptic and bielliptic. For example \cite{KuKo}, $C=\{y^2=(x^2-1)(x^2-\lambda^2)(x^2-\mu^2)(x^2-\lambda^2/\mu^2)\}
,\,\,\,\lambda^2 \neq \mu^2,\,\,\lambda^2,\,\mu^2 \in \C\setminus \{0,\,1\},$ 
is such a curve. 
\end{rem}

\medskip

\begin{ex} A genus 6 nodal plane sextic $C$ is neither hyperelliptic, nor bielliptic \cite[V.A-12]{ACGH}. 
Hence by Proposition \ref{hypersym3}, $V'=C_2$ is a hyperbolic surface. 

As for examples of plane quartics with the same property, a
genus 3 curve $C$ is hyperelliptic or bielliptic iff it possesses a
holomorphic involution, as noted in Remark \ref{ell}(c) above. 
Or equivalently, iff the automorphism group Aut$\,C$
contains an element of order 2, 
that is, iff Aut$\,C$ is of even order. 
In fact, in \cite{KuKo} all 
the genus 3 curves were classified according 
to the order of the group Aut$\,C$. 
E.g., this order equals 9 for the 
plane quartic $C=\{x^4-xz^3-y ^3z=0\}$ and whence, $C$ 
is neither hyperelliptic nor bielliptic.

\end{ex}
 
\medskip Next we show that the symmetric power $C_d,\,\,d\ge 2,$ of a curve
$C$ is never C-hyperbolic and indeed is just the opposite in the sense
of the following definition.

\medskip

\begin{defin} A complex space $X$ is called {\it Liouville} 
if it carries no nonconstant bounded holomorphic functions
\cite{Li}. We call it {\it super-Liouville} if 
the universal cover over $X$ is Liouville. \end{defin}

\medskip

Evidently, any compact complex space is Liouville; 
in fact, this property is opposite 
to being Carath\'eodory hyperbolic. 
In turn, being super-Liouville is 
opposite of being C-hyperbolic.
By Lin's Theorem \cite{Li}, 
any compact complex space
with nilpotent (or nilpotent-by-finite) fundamental group
is super-Liouville. 

In particular, this is so for complex tori, as well 
as for manifolds birational to complex tori. Thus, any 
cover over the symmetric power $C_g$ of a 
genus $g$ curve $C$ is Liouville 
(indeed, the Abel-Jacobi map $u_g\,:\,C_g \to J_C$ is birational). 

Another example of a super-Liouville variety is provided by the theta-divisor
$\Theta \subset J_C$ of the Jacobian of a genus $g\ge 3$ curve $C$  with
general moduli. Indeed, the divisor $\Theta$ being ample (see
\cite{Mu}), by the Lefschetz hyperplane section theorem (see for
instance, \cite{Mi}), the embedding $\Theta \hookrightarrow J_C$ induces
an isomorphism of the fundamental groups.
As for the symmetric powers, we have the following statement:

\medskip \begin{lem} For any curve $C$ of genus $g$ and for any 
$d\ge 2$ we have
$\pi_1(C_d) \simeq \Z^{2g}.$ \end{lem}

\medskip 

\noindent \proof It is known that $H_1(C_d;\,\Z) \simeq \Z^{2g}$
\cite{McD}. Therefore, it suffices to show that the fundamental group
$\pi_1(C_d)$ is abelian. Let $\pi\,:\,C\times \dots \times C\to C_d$
be the natural projection. We first show:

\medskip
\noindent {\it Claim: $\pi_*\,:\,\pi_1(C \times \dots \times C) 
\to \pi_1(C_d)$ is a surjection.}

\medskip We denote by $\Delta=\bigcup_{i\neq j} \Delta_{ij}$ the union of the diagonal hypersurfaces in 
the cartesian power $C^d=C\times \dots \times C,$ and we let 
$\Delta' = \pi(\Delta),\,\,\,$ 
$\pi^0 = \pi\,|\,(C^d \setminus \Delta)\,:\,C^d\setminus \Delta \to  C_d
\setminus \Delta'.$ Since $\pi^0$ is a non-ramified covering, we have the
exact sequence  $${\bf 1}\longrightarrow \pi_1(C^d \setminus \Delta)
\stackrel{\pi^0_*} {\longrightarrow} 
\pi_1(C_d \setminus \Delta') \stackrel{\rho} 
{\longrightarrow} S_d {\longrightarrow} {\bf 1}$$ 
where $S_d$ stands for the symmetric group.\footnote
{The group $\pi_1(C_d \setminus \Delta')=B_{g,d}$ is called the 
$d$-th braid group of a genus $g$ compact Riemann surface.}
This sequence can be   extended 
to the commutative diagram

\hspace{1.3cm}\begin{picture}(200,210)

\put(135,180){$\pi_1(C^d \setminus \Delta) \quad\stackrel{i_*}
{\longrightarrow} \,\,\,\pi_1(C^d) \,\longrightarrow {\bf 1}$}

\put(154,165){$\vector(0,-1){25}$}
\put(160,150){$\pi^0_*$}
\put(250,165){$\vector(0,-1){25}$}
\put(255,150){$\pi_*$}

\put(40,120){${\bf 1}\longrightarrow N' \quad\longrightarrow \quad 
\pi_1(C_d \setminus \Delta') \quad \stackrel{i'_*}{\longrightarrow}
\,\,\,\, \pi_1(C_d)\longrightarrow {\bf 1}$}

\put(78,105){$\vector(0,-1){25}$}
\put(83,90){$\rho|N'$}
\put(154,105){$\vector(0,-1){25}$}
\put(160,90){$\rho$}
\put(250,105){$\vector(0,-1){25}$}

\put(40,60){${\bf 1}\longrightarrow S_d \quad\stackrel{\rm
id}{\longrightarrow} \qquad \,\,S_d \,\,\qquad
\longrightarrow \qquad\,\,\, {\bf 1}$}

\put(78,45){$\vector(0,-1){25}$}
\put(154,45){$\vector(0,-1){25}$}

\put(75,0){${\bf 1} $}
\put(151,0){${\bf 1}$}
\end{picture}

\bigskip\bigskip \noindent where  $i\,:\,C^d \setminus \Delta
\hookrightarrow C^d$, $i'\,:\,C_d \setminus \Delta'\hookrightarrow C_d$  are
the natural embeddings and $N' = {\rm ker}\,i'_*$.
This kernel $N'$  is generated, as a normal subgroup, by a
vanishing loop $\alpha'$ of the diagonal
$\Delta'$  in the manifold $C_d$.  We note that $\rho(\alpha')$ is a
transposition in $S_d$.  Indeed, the kernel of $i_*$ is generated, as a
normal subgroup, by vanishing loops $\alpha_{ij}$ of the diagonals
$\D_{ij}$ in $C^d$; if we choose $\alpha'$ so that, for instance,
$\pi^0_*(\alpha_{12}) = (\alpha')^{2}$, then $\rho(\alpha')=(1\,2)$. Since
$\rho(N')$ is a normal subgroup of $S_d$ containing a transposition, we
conclude that $\rho(N')=S_d\,$; i.e., the first two columns, as well
as the rows, of the diagram are exact.  By the usual diagram chasing (as in
the ``nine lemma"), we conclude that the last column is exact, too; i.e.,
$\pi_*$ is surjective. This proves the claim. 

\smallskip
Now we can prove that the group $\pi_1(C_d)$ is abelian.  Indeed, let
$\{a_i^{(k)}\}_{1\le i\le 2g}$ denote the set of standard generators of the
$k$-th factor of the fundamental
group  $\pi_1(C^d) \simeq \pi_1(C) \times\dots \times \pi_1(C)$; we
have  $[a_i^{(k)},\,a_j^{(l)}]=1$ for all $i,j=1,\dots,2g$, if $k\neq l.$ 
Hence, since $d\ge 2$ and
$\pi_*(a_i^{(k)})=\pi_*(a_i^{(l)})=:a_i'\in \pi_1(C_d),
\,\,\,i=1,\dots,2g,$ we have
$1=[\pi_*(a_i^{(k)}),\,\pi_*(a_j^{(l)})] = [a_i',\,a_j']$ for all
$i,j=1,\dots,2g.$ But by the above claim, 
$$\{a_i'=\pi_*(a_i^{(k)})\}_{1\le i\le 2g}$$ 
is a set of generators of
the group $\pi_1(C_d)$. Thus, this group is, indeed, abelian. \qed

\medskip

\begin{cor} \label{super-L} For any curve $C$ and any $d\ge 2$ 
the symmetric power $C_d$ of $C$
is super-Liouville,
that is, any non-ramified covering over $C_d$ is Liouville. 
\end{cor}

\medskip Next we show how to produce examples of hyperbolic surfaces in
$\PP^3$ starting with the symmetric square of a generic genus $g \ge 3$ curve
$C$.  To obtain these examples, we need to know the relations between 
the
divisor theories on the cartesian square $V=C\times C$
and on the symmetric square
$V'=C_2.$  We let $\pi:V\to V'$ denote the projection.
If $D'$ is a divisor on $V'$, it follows from the push-pull formula that
\begin{equation}\label{pi^*}  \pi_*\pi^* D' = 2D'\qquad 
{\rm and}\qquad \pi^*D'\cdot A=D'\cdot \pi_*A\,,
\end{equation} for any divisor $A$ on $V$.
Thus, if $D_1,\,D_2$ are divisors on $V',$ we have
\begin{equation}\label{pi_*}
\pi^*D_1\cdot \pi^*D_2 = 2 D_1\cdot
D_2\,.\end{equation}

We let $E',\ \D'$ be the divisors on $V'$ given by $\pi^*E'=E+F$,
$\pi^*\D'=2\D.$ It is well known that (when $C$ has general moduli) the
Neron-Severi group $\NS(V')$ is generated by the classes of $E'$ and $\Theta'$
where $\Theta'$ is the pull-back of the theta divisor $\Theta$ on $J_C$ by the
Abel-Jacobi map $u_2$ (see \cite{ACGH, GriHarr}).\footnote{Recall that by
Lefschetz' Theorem (see \cite{Pir} for a modern proof), for a curve $C$ with
general moduli, the Neron-Severi group $\NS(J_C)$ is a free abelian group
generated by the class of the theta divisor $\Theta$.} Furthermore
\cite{McD} (see also \cite[VIII.(5.4)]{ACGH}, \cite{Kou}),
\begin{equation}\label{deltatheta}\D'\equiv (2g+2)E'-
2\Theta'\,,\end{equation} and hence    $E',\ \half \D'$   generate
$\NS(V')$,  where  $\half \D'\equiv (g+1)E'-\Theta'.$ We
write 
\begin{equation} \label{D'} D'_{a,k}=(a+k)E'-{k\over 2}\D'\,,
\end{equation}
so that    
$E'\equiv D'_{1,0},\,\,\,\Delta'\equiv 2D'_{1,-1},\,\,\,\Theta'\equiv 
D'_{g,1}$, and also
\begin{equation} \label{*'} \pi^*D'_{a,k}=D_{a,k}\,.\end{equation}
By (\ref{pi_*}), we have 
\begin{equation}\label{proD'} D'_{a,k}\cdot D'_{c,l} = ac -gkl\,.
\end{equation}
We now compute $K:=K_{V'}.$    Since $\Delta$ 
is the ramification divisor of the branched 
cover $\pi\,:\,C\times C \to C_2$ we have   
$$\pi^*K=K_V-\D\equiv D_{2g-3,1}\,,$$
and so by (\ref{D'}) (cf.\ \cite[VIII.(5.4)]{ACGH}),
\begin{equation} \label{K'} K \equiv D'_{2g-3,1}\equiv (2g-2)E'-\half\Delta'
\,.\end{equation} 
As for the topological invariants of the surface
$V',$ we have: \begin{eqnarray} \label{c1'} c_1^2=K^2 &=&
(2g-3)^2-g=4g^2-13g+9\\ c_2=e(V') &=& \half[e(V)+e(\D)]\nonumber\\
\label{c2'} &=& \half[e(C)^2+e(C)]\nonumber\\
&=& (2g-3)(g-1)=2g^2-5g+3\,.\end{eqnarray}
Hence, 
\begin{equation} \label{chi'} \chi({\cal O}_{V'})={c_1^2+c_2 \over 12} 
= {(g-1)(g-2)\over 2}\,.\end{equation}

Next we describe some of the very ample divisors on the surface $V'$.
For the projection $\pi:V\to V'$ we write $\pi(z,w)=\{z,w\}.$ If $A$ is a
divisor on $C$, we let $A^{(2)}$ be the (unique) divisor on $V'$ given by
$$\pi^*A^{(2)}=A\times C + C\times A\,.$$ If $\wt A$ 
is linearly equivalent to
$A$, then $\wt A^{(2)}$ is linearly equivalent to $A^{(2)}$; thus for a
holomorphic line bundle $L\to C$, we have a unique line bundle $L^{(2)}\to
V'$, and $\pi^*L^{(2)}=L{\widehat\otimes}L:= {\rm pr}_1^* L + {\rm pr}_2^* L.$

It follows from  the
Nakai-Moishezon criterion and (\ref{pi_*}), (\ref{pi^*}) that $D$ is ample iff
$\pi^*D$ is.  Thus if $L$ is a line bundle on $C,$ then $L$ is ample iff
$L{\widehat\otimes}L$ is ample iff $L^{(2)}$ is ample. We further note that
the latter holds also for ``very ample":

\begin{lem}\label{vaL2} $L\to C$ is very ample iff $L^{(2)}\to V'$ is very
ample.\end{lem}

\noindent \proof Let $\iota:C\to V'$ be given by $\iota(z)=\{z,x_0\}.$  If  $L^{(2)}$
is very ample, then so is $L=\iota^*L^{(2)}.$ Conversely, suppose that $L$
is very ample. We first show that $L^{(2)}$ separates points. For $s,t\in
H^0(C,L)$, we let $s\odot t$ denote the section of $L^{(2)}$ given by
 $$s\odot t(\{z,w\})=s(z)t(w)+s(w)t(z)\,.$$  
Let $\zeta\ne\eta\in V'$ be two
arbitrary points. We must find a section $\lambda\in H^0(V',L^{(2)})$ such that
$\lambda(\zeta)=0,\ \lambda(\eta)\ne 0.$ First suppose that
$\zeta=\{z_1,w_1\},\ z_1\not\in\{z_2,w_2\}=\eta.$  Since $L$ is very ample,
we can find a section $s\in H^0(C,L)$ such that $s(z_1)=0,\ s(z_2)\ne 0,\
s(w_2)\ne 0.$  Then $\lambda=s\odot s$ is our desired section.  The other
possibility is that $\zeta=\{z_1,z_1\},\ \eta=\{z_1,z_2\}.$ In this case, we
let $\lambda=s_1\odot s_2$ where $s_i(z_j)$ vanishes if $i=j$ and is
nonzero if $i\ne j.$ Thus, $L^{(2)}$ separates points.

Next consider an arbitrary nonzero tangent vector $X\in T_\zeta(V'),\
\zeta=\{z_0,w_0\}\in V'.$  To complete the proof, we must find a section
$\lambda\in H^0(V',L^{(2)})$ such that $\lambda(\zeta)=0$ but the 1-jet 
${\cal J}_1\lambda$ of $\lambda$ does not vanish in the $X$ direction. First
assume that $z_0\ne w_0$ and write $X=c_1{\partial\over\partial
z}|_{z_0} + c_2 {\partial\over\partial w}|_{w_0}$; we may assume that
$c_1\ne 0.$ We can then let $\lambda=s\odot s$, where $s$ is a section in
$H^0(C,L)$ such that $s(z_0)=0$, $s(w_0)\ne 0$, and ${\cal J}_1 s|_{z_0} \ne
0.$ Now consider the case $\zeta=\{z_0,z_0\}.$  Choose a local frame
$e$ for $L$ at $z_0$, and use a local coordinate centered at $z_0$ so that
we may write $z_0=0.$ Since $(z+w,zw)$ are local coordinates at $\zeta\in
V'$, it suffices to find $\lambda\in
H^0(V',L^{(2)})$ such that \begin{equation}\label{j1}\lambda=
(0+c_1(z+w)+c_2 zw +\cdots)e\odot e\end{equation} for arbitrary
$c_1,c_2\in\C.$ Choose $s,t\in H^0(C,L)$ with ${\cal J}_1s=ze$, ${\cal J}_1
t=(c_1+{c_2\over 2}z)e.$ Then the section $s\odot t\in H^0(V',L^{(2)})$
satisfies (\ref{j1}). Thus $L^{(2)}$ is very ample. \qed

\medskip

\begin{cor} \label{newcor}
Let $C$ be a non-hyperelliptic curve of genus $g\ge 3.$ Then the divisor
$K_C^{(2)}\equiv D'_{2g-2, 0}$ on the symmetric square $V'=C_2$ of $C$ is very
ample.
\end{cor}

\noindent \proof Indeed, since the canonical divisor $K_C$ is very ample, 
by  \ref{vaL2}, $K_C^{(2)}$ is very ample, too. \qed

\medskip
\begin{prop}
For a very ample divisor $H\equiv D'_{a,k}$ on the surface $V'$, 
in the notation from Section~1, we have:
\begin{eqnarray}\label{nsym} n &=& a^2-gk^2\\
\label{gHsym}  g_H &=& \half[n+2 +(2g-3)a-gk]\\
\label{gSsym} g_\bS &=& \half (n^2-7n-7g^2+25g+8) +(n-12) g_H\\
\label{bsym} b &=&{1\over 2}n(n-5) - 5(g_H-1)- (g-1)(g-3)
\end{eqnarray}\end{prop}

\noindent \proof 
These follow as before from  (\ref{b1}), (\ref{adj}), (\ref{g}) 
and the above formulas (\ref{D'})--(\ref{chi'}).  For example, 
(\ref{gHsym}) follows from  the adjunction formula $2g_H-2 =H\cdot K +n$,
where
$$H\cdot K = D'_{a,k}\cdot D'_{2g-3,1} =
(2g-3)a-gk\,.$$ \qed

\begin{thm}\label{deg16} Let $C$ be a genus $3$ curve that is neither
hyperelliptic nor bielliptic, and consider the divisor $H=K_C^{(2)}$ on
$V'=C_2$.  Then $H$ is very ample and a generic projection $\bV'\subset
\PP^3$ of the image of the projective embedding $V'\hookrightarrow \PP^N$
given by $H$ is a Kobayashi hyperbolic surface of degree $16$ in $\PP^3.$
\end{thm}

\noindent \proof By Corollary \ref{newcor}, the divisor $K_C^{(2)}\equiv
D'_{4,0}$ is very ample. Furthermore by (\ref{nsym})--(\ref{gSsym}), we have
$n=16$, $g_H=15,\,\, g_{\bS}=142$.  Now the conclusion follows from
Corollary \ref{khyp} and Proposition \ref{hypersym3}. \qed

\medskip 

\begin{rem} By \cite{Za 1}, 
to obtain a {\it smooth} hyperbolic surface in $\PP^3$ of degree 16, it 
is enough to perturb a little the coefficients of 
the equation which defines 
the singular hyperbolic surface provided by Theorem \ref{deg16}. 
\end{rem} 

\medskip The following proposition tells us that actually, $16$ is the lowest
possible degree of a projective embedding of the symmetric square
$V'=C_2$ of a generic genus $g\ge 3$ curve $C.$

\smallskip

\begin{prop} \label{estsym} Let $C$ be a genus $g\ge 3$ curve with general
moduli. Then $\delta(C_2)=16$ for
$g=3$, $\delta(C_2)\ge 20$ for $g=4$, and  \begin{equation}\label{insym}
\delta(C_2) > \sqrt{2g(g+11)}+5\quad \mbox{for}\quad g \ge
5\,.\end{equation} \end{prop} 

\noindent \proof Let $H\equiv D'_{a,k}$ be a hyperplane section of
$V'=C_2\hookrightarrow \PP^N.$ By (\ref{nsym}), $|k| < a/\sqrt{g}.$ Since
$a=\deg\iota^*H\ge \delta(C)$, where $\delta(C)$ is as in
Proposition~\ref{est}, it then follows from (\ref{gHsym}) and 
(\ref{delta}) that
\begin{equation} \label{gHinsym} g_H > \half \left[n+2
+\left(2g-\sqrt{g}-3\right)\delta(C)\right]\,.\end{equation} 
On the other hand, 
from
(\ref{bsym}) and the fact that $b$ is nonnegative, we have the inequality
\begin{equation}\label{in2sym} g_H \le {n^2 - 5n + 10 - 2(g-1)(g-3)\over
10}\,. \end{equation} Combining (\ref{gHinsym}) and (\ref{in2sym}), we
obtain  
\begin{equation}\label{insym2}
(n-5)^2> 2g^2-8g+31 +5\delta(C)(2g-\sqrt{g} -3)\,.
\end{equation}

We begin with the case $g=3$. Then $\delta(C)=4$, and
(\ref{insym2}) yields $(n-5)^2 > 50$ or $n\ge 13$. But the only solutions of
the diophantine equation $3k^2=a^2-n$ (provided by (\ref{nsym})) with
$n=13,14,15$ and $0<a\le 11$ are $n=13$, $(a,k)=(4,1),(5,2)$ and
$(11,6)$.  In the first two cases $g_H=12$ and in the third case $g_H=15$.
On the other hand, (\ref{in2sym}) with $n=13,\ g=3$ yields $g_H\le 11$.
Therefore, if $n\le 15$,  we can replace $\delta(C)$ with
12 in (\ref{insym2}), which yields 
$$100\ge (n-5)^2 > 25+ 60 (3-\sqrt{3}) > 100\,,$$ 
a contradiction.  Hence $\delta(C\times C)\ge 16$. However,
Theorem~\ref{deg16} implies that $\delta(C\times C)\le 16.$ 

If $g\ge 4$, then $\delta(C)\ge 5$ and (\ref{insym2}) yields 
\begin{equation}\label{insym3}  (n-5)^2 >
2g^2+42g-25\sqrt{g}-44 \,.\end{equation}
For $g=4$, (\ref{insym3}) gives
$(n-5)^2 >106$
or $n\ge 16$. But for $g=4$, equation (\ref{nsym}) becomes
$$(a-2k)(a+2k)=n\,,$$ which has no integer solutions for $16\le n\le
19$.  Thus $n\ge 20.$

For $g\ge 5$, (\ref{insym3}) and the inequality
$25\sqrt{g}+44 < 20g$ yields (\ref{insym}). \qed

\medskip

\begin{cor} \label{cinsym} Let $C$ be a genus $g\ge 3$ curve with general
moduli. Then for any projective embedding
$C_2\hookrightarrow \PP^N$, we have $g_\bS\ge 130$.\end{cor}

\noindent \proof  For the case $g=3$, using
(\ref{gHinsym}) with $n\ge 16$ and $\delta(C)=4$, we get  $g_H
>15-2\sqrt{3}>11.5$ (i.e., $g_H\ge 12$)  and thus by (\ref{gSsym}),
$g_\bS\ge 130$. Similarly for $g=4$, using $n\ge 20$ and
$\delta(C)\ge 5$, we get $g_H\ge 19$ and $g_\bS\ge 280$.

We now assume that $g\ge 5$. {From} (\ref{insym}), we obtain
\begin{equation}\label{insym4} n > \sqrt{2}\, g +10 \,.\end{equation} 
The inequalities (\ref{gHinsym}), (\ref{delta}) and (\ref{insym4}) yield 
\begin{equation}\label{ghbound} g_H > 
\sqrt{2}\,g^{3/2}+{3\over 2}g- {3\over 4}(1+2\sqrt{2}) 
g^{1/2}+{15\over 4} >
\sqrt{2}\,g^{3/2}+4\,.\end{equation} It then follows from (\ref{gSsym}), 
(\ref{insym4}),
and (\ref{ghbound}) that
$$g_\bS > 2 g^{5/2}-{5\over 2}g^2-2\sqrt{2}\,g^{3/2}+{21\sqrt{2}+25 \over
2}g+11\ge  165\,.$$ \qed

\medskip Finally, from Corollaries \ref{khyp}, \ref{cinsym} and
Propositions
\ref{hypersym3}, \ref{estsym}, we obtain the following theorem.

\medskip 
\begin{thm} \label{hypsym} Let $V' = C_2$ where $C$ 
is a smooth projective curve of genus $g \ge 3$ with general moduli.
Then a generic projection $\bV' \subset \PP^3$ of 
the image of any projective embedding $V' \hookrightarrow \PP^N$ 
is
a Kobayashi hyperbolic surface of degree $n \ge 16.$\end{thm}

As an aside, we have the following observation 
related to Lemma \ref{gamma} and Remark \ref{ell} above:

\begin{prop}  \label{gammasq} 
Let $f:C\to R$ be a $2:1$ morphism from a curve $C$ of genus $g\ge
2$ to another curve $R$, and let
$\G=f^{\vee} (R) \subset C_2$, where $f^{\vee}$ is as in Remark
$\ref{ell}$ above.  Then $\G^2=1-g$, and if $R$ is
nonrational,  then the Picard number of the surface $C_2$ is at
least $3$. This is the case, in particular, if $C$ is bielliptic.
\end{prop}

\medskip

\noindent\proof It is clear that $\G$ is a smooth curve isomorphic to $R$.
Thus by the adjunction formula and (\ref{K'}), 
\begin{equation}\label{f2} \G^2=2g_R-2 - K\cdot \G=
2g_R-2 -(2g-2)\G\cdot E' -\half\G\cdot \Delta'\,,\end{equation}
where $g_R$ denotes the genus of $R$. Here $\G\cdot E'=1$ and
by the Riemann-Hurwitz formula,  
\begin{equation}\label{f3} \Gamma\cdot\Delta'=2+2g-4g_R\,.\end{equation}
Hence by (\ref{f2}) and (\ref{f3}), $\Gamma^2 =1-g$.

To obtain the second statement, suppose on the contrary that
$E'$, $\half\D'$ generate $\NS(V')$, and write $\G=D'_{a,k}$.  Then
$a=\G\cdot E'=1$ and $$2+2g-4g_R=\G\cdot\D'=D'_{1,k}\cdot D'_{2,-2}=
2+2gk\,,$$ so that $k=(g-2g_R)/g$.  Then
$1-g=\G^2=1-(g-2g_R)^2 /g$, which gives $g=g_R$ (recall that $g_R >
0$). But then by (\ref{f3}), $g\ge 2g_R-1=2g-1$ or $g\le 1$,
contradicting our assumption.
\qed 

\medskip

\noindent\begin{rem} Let $C$ be a genus 3 curve. It is well known 
(see \cite[VI.\S 4]{ACGH}) that the theta-divisor 
$\Theta \subset J_C$ is singular iff $C$ is hyperelliptic. 
According to Lemma \ref{gamma} and Proposition \ref{gammasq}, 
if $C$ is hyperelliptic then there is a unique (smooth) 
rational curve $\G$ on the surface $V'=C_2$, and $\G^2=-2$,
so that the Abel-Jacobi morphism 
$$u_2\,:\,C_2 \to W_2 \simeq \Theta \subset J_C$$
contracts $\G$ to a (unique) singular point of type $A_1$ of
the surface $W_2 = \Theta + p,\,\,p\in J_C.$ \end{rem}

\bigskip 
\section{Divisors on cartesian and symmetric squares of curves of genus
$\ge 2$}\label{sec-divisors}

To provide more examples of the hyperbolic surfaces in $\PP^3$ given by
Theorems \ref{hyp} and \ref{hypsym}, we must find projective embeddings of
$V$, resp.\ $V'$, i.e., we must find sufficient conditions for a divisor
$H\equiv D_{a,a',k}$, resp.\ $H\equiv D'_{a,k}$, to be very ample. In
this section we use a description of the nef cones on symmetric
squares given in Kouvidakis \cite{Kou}
together with Reider's
characterization of very ampleness \cite{Rei} to show that $H$ is very ample
in the following cases:

\begin{enumerate}
\item[i)] genus\;$C=2$, $H\equiv D_{a,k}\in \NS(C\times C),\ a\ge
5$ and $2|k|\le a-3$ (Theorem \ref{g2ample}(e)).

\item[ii)] $C$ is a generic curve of genus 3, $H\equiv D_{a,k}\in \NS(C\times
C),\ a\ge 7,\ (a,k)\ne (8,2)$ and $-{1\over 3} (a-4) \le k \le {5\over
9}(a-4)$ (Theorem \ref{veryample3}).

\item[iii)] $C$ is a generic curve of genus 3, $H\equiv D'_{a,k}\in \NS(C_2),\
a \ge 7$, $(a,k)\ne (7,3)$ or
$(9,4)$, and $4-a \le 3(k-1) \le {1\over 3}(5a-16)$
(Theorem \ref{g2ample'}(b$'$)). 

\item[iv)] $C$ is a generic curve of genus 4, $a\ge 9$ and $9-a \le 4k \le
2a-10$ (Theorem \ref{g2ample'}(b$''$)). 
\end{enumerate}

We let $\Sigma'$ denote the subgroup of the 
Neron-Severi group $\NS(V')$
generated by the classes of 
$E',\;\half\Delta'$; for a curve with general moduli, $\Sigma'=\NS(V')$.
We begin by restating (using the basis $\{D'_{1,0}=E',\ 
D'_{0,1}\equiv \Theta'-gE'\}$ of $\Sigma'\otimes\Q$) Kouvidakis's
description \cite[Thm. 2]{Kou} of the effective and 
nef cones on the symmetric square of a generic curve:

\begin{thm}{\bf (Kouvidakis)} \label{Kou}{\it Let $V'=C_2$ 
be the symmetric square of a genus $g\ge 2$ curve with general moduli.
Denote by ${\rm EFF}(V')$, resp.\ ${\rm NEF}(V')$, the cone of 
(the classes of)
quasi-effective,
resp.\ nef, $\Q$-divisors in 
$\Sigma' \otimes \Q$. 
We also set
\begin{displaymath}{\cal E} = \left\{
D'_{a,k} \in \Sigma'\otimes \Q\,\left\vert
\begin{array}{ccc}
-a&\le& k\le a \qquad \,\,\,\,\,\,\,\,\,\,\,{\rm if} \qquad g=2\\
-a&\le& k\le 3a/5 \qquad \,\,{\rm if} \qquad g=3\\
-a&\le& k\le a/\sqrt{g} \qquad {\rm if} \qquad g\ge 4
\end{array}\right\}\right .
\end{displaymath}
and $$ {\cal E}' = \left\{D'_{a,k} \in \Sigma'\otimes \Q\,\left\vert\,
-a\le k \le {a\over \sqrt{g}-1}\right.\right\}$$
for $g\ge 5;$ 
further, 
\begin{displaymath}{\cal N} = 
\left\{
D'_{a,k} \in \Sigma'\otimes \Q\,\left\vert
\begin{array}{ccc}
-a/2&\le& k\le a/2 \qquad \,\,\,\,\,{\rm if} \qquad g=2\\
-a/3&\le& k\le 5a/9 \qquad\,\, {\rm if} \qquad g=3\\
-a/g&\le& k\le a/\sqrt{g} \qquad {\rm if} \qquad g\ge 4
\end{array}
\right\}\right .
\end{displaymath}
and $$ {\cal N}' = \left\{D'_{a,k} \in \Sigma'\otimes \Q\,\left\vert\,
-{a\over g}\le k \le a{\sqrt{g}-1\over  g}\right.\right\}\,$$
for $g\ge 5.$ Then we have:
$${\cal E} \subseteq {\rm EFF}(V')  \subseteq {\cal E}'\,,\qquad {\rm resp.} 
\qquad {\cal N}' \subseteq {\rm NEF}(V')  \subseteq {\cal N}\,,$$
and $${\rm EFF}(V') = {\cal E}\,,\qquad {\rm resp.} 
\qquad {\rm NEF}(V') = {\cal N}\,,$$
if $g=2,\,3$ or if $\sqrt{g} \in \Z.$ }\end{thm}

\medskip
Here we call a divisor $D$ {\it quasi-effective\/} iff $mD\equiv G$ for
some $m>0$ and for some effective divisor $G$. Recall that the ample
cone (sometimes called the {\it K\"ahler cone\/}) is the 
interior of the nef cone, and hence it can be
described by making the inequalities in the definition of ${\cal
N},\,{\cal N'}$ strict. 

We let $\Sigma^\sym=\pi^*\Sigma'$ denote the
subgroup of $\Sigma\subset \NS(V)$ generated by the ``symmetric" divisors
$\{D_{a,k}\}$. The following elementary observation allows us to transfer
Kouvidakis's description of the effective and nef (as well as ample) cones
to the case of symmetric divisors $D_{a,k}\in\Sigma^\sym\subset \NS(V)$.

\begin{prop}\label{elementary} Let $V$, resp.\ $V'$ 
be the cartesian, resp.\
symmetric, square of a smooth projective curve $C$.  
Then the class $D_{a,k}\in \Sigma^\sym \subset NS(V)$ is
quasi-effective, resp.\ nef, resp.\ ample, iff
$D'_{a,k}\in\Sigma'\subset NS(V')$ is likewise. \end{prop}

\noindent \proof The proposition is an immediate consequence of the
push-pull formulas (\ref{pi^*}), (\ref{*'}) and the Nakai-Moishezon criterion.
\qed 

\medskip Let $C$ be a smooth projective curve of genus $g = 2$.
Then $C$ is a hyperelliptic curve, and there is a unique
hyperelliptic involution $i\,:\,C \to C.$ The 
orbits $(p,\,i(p))$ of the involution $i$ 
are divisors of the canonical 
class, and the six
fixed points of $i$ are 
the Weierstrass points of the curve $C.$ 
Let $\varphi=\varphi_{|K_C|}\,:\,C \to \PP^1$ be 
the quotient by the involution $i.$ 
The involution $i$ or, equally, the morphism $\varphi$
defines the symmetric correspondence (the graph of $i$)
$$D(i) \subset V = C \times C,\,\,
D(i) := \{(p,\,i(p))\,|\,p \in C\}$$ 
of valence $1$ and of bidegree $(1,\,1).$ 
Hence, $$D(i) \equiv D_{1,1} = 2E + 2F - \Delta \in \Sigma\,,$$
whereas $K = K_V \equiv D_{2,0}=2E + 2F.$ 
We have $\Delta^2 = D(i)^2 = -2.$ 

To obtain a description of the very ample (as well as some of the
globally generated) divisors of type $D_{a,k}$ on the surface $V$, 
resp.\ $V'$, we shall
use  Reider's characterization of global generation and of very ampleness
\cite[Theorem 1(ii)]{Rei} (see also \cite{Laz}):

\medskip

\begin{thm}{\bf (Reider)} \label{re}{\it Let $L$ be a nef line bundle 
on a smooth projective
surface $V$ such that $L^2 \ge 5$, resp.\ $L^2 \ge 10.$ 
Then the adjoint line bundle
$K+L$ is globally generated, resp.\ very ample, 
unless there exists an effective divisor $\G$ 
on $V$ which verifies one of the following conditions $(i)$--$(ii)$, resp.\ 
$(i')$--$(iii')$:

\noindent $(i)$ $L\cdot \G = 0 $ and $\G^2 = -1 ;$ \hfill 
$(i')$ $L\cdot \G = 0 $ and $\G^2 = -1 $ or $-2;$ 

\noindent $(ii)$ $L\cdot \G = 1 $ and $\G^2 = 0.$ \hfill
$(ii')$ $L\cdot \G = 1 $ and $\G^2 = -1 $ or $0;\,\,\,$ 

\noindent \hfill $(iii')$ $L\cdot \G = 2 $ and $\G^2 =
0.\,\,\,\,\,\,\,\,\,\,\,\,\,\,\,\,\,\,\,$ } \end{thm}

In Theorem \ref{g2ample} below, we describe all the quasi-effective, nef, and
very ample, as well as some globally generated, divisors of type $D_{a,k}$ on
the surface $V=C\times C$, for an {\it arbitrary} genus 2 curve $C$. For a
{\it generic} curve $C$, statements (a)--(c) of Theorem \ref{g2ample} and
Corollary \ref{g2'ample} can be obtained (in view of Proposition
\ref{elementary}) from the genus 2 case of Kouvidakis's Theorem
\ref{Kou}, but we provide a direct proof for the reader's convenience.

\begin{thm} \label{g2ample} Let $V = C \times C$ 
be the cartesian square of a genus $2$ curve $C$. 
Then a divisor $H \equiv D_{a,k} \in \Sigma^\sym\subset\NS(V)$
is

\noindent {\rm (a)} quasi-effective iff 
\begin{equation} \label{effective} |k| \le a\,;\end{equation}
{\rm (b)} nef iff 
\begin{equation}\label{nef} 2|k| \le a \,;\end{equation}
under this condition it is also big unless $H \equiv 0;$

\smallskip\noindent {\rm (c)} ample iff \begin{equation}\label{a} 2|k| \le
a-1\,;\end{equation} 
{\rm (d)} ample, globally generated, and non-special if $(a,k)=(4,0)$ or $a
\ge 5$ and  \begin{equation}\label{gg} 2|k| \le a-2\,;\end{equation}

\smallskip\noindent {\rm (e)} very ample iff $a \ge 5$ and 
\begin{equation}\label{va} 2|k| \le a-3\,.\end{equation}\end{thm}

\noindent \proof For a divisor class $\G\in \NS(V)$, we write $\G=\G^\Sigma+
\G^\perp$, where $\G^\Sigma\in \Sigma\otimes\Q$ and $\G^\perp\in 
\NS(V)\otimes\Q$ is
such that $ \G^\perp\cdot D=0 \ \forall D\in\Sigma$.  We first establish the
following

\smallskip\noindent {\it Claim:\/} Suppose $\G\in\NS(V)$ is effective
and non-zero, and let $\G^\Sigma\equiv D_{c,c',l}$. Then 

\smallskip 

\noindent (i) $c,c',4l\in\Z$,
$c,c'\ge 0$ and $c+c'\ge 1$; 

\smallskip 

\noindent (ii)  if $\G \neq \Delta$ or $D(i)$ then
\begin{equation}\label{qeff2} 4|l|\le
c+c'\,;\end{equation} 
(iii) furthermore, $c+c'\ge 2$ if $\G\not\equiv E$ or $F$.

To verify the claim, we can assume without loss of generality that $\G$ is an
irreducible curve different from $\D \equiv D_{1,-1},\ D(i) \equiv
D_{1,1},\ \C\times\{p\}\equiv E$ and $\{p\}\times C \equiv F$.  
First, we note that $\G\cdot E
=\G^\Sigma\cdot E=c' \in\Z$. Since $\G\not\equiv E$, $c' \ge 1$;
indeed, if $c'=0$, then $\G\cap E=\emptyset$, which can only happen if $\G=
C\times\{p\}$, contrary to our assumption.
Likewise $\G\cdot F=c \in \Z^+$, and furthermore, $\G\cdot
D_{0,-1}=4l\in\Z$. To show (\ref{qeff2}), we note that $$0\le \G\cdot
D_{1,\pm 1} = D_{c,c',l}\cdot D_{1,\pm 1}=c+c' \mp 4l\,,$$ 
which gives (\ref{qeff2}). Thus the
claim is established.

\smallskip
(a): If (\ref{effective}) holds, then we easily see that
$D_{a,k}$ is an effective linear combination of $E+F,\ \D$ and $D(i)$.
The converse is an immediate consequence of the above claim.

\smallskip
(b): If $H$ is nef, then $$0\le H\cdot D_{1,\pm 1} = 2(a\mp 2k)$$ and hence
(\ref{nef}) holds. Conversely, assume (\ref{nef}), and let $\G\in\NS(V)$ be an
effective divisor. Write $\G^\Sigma=D_{c,c',l}$; by the claim,
\begin{equation}\label{HGamma}  H\cdot\G= H\cdot \G^{\Sigma} = H\cdot D_{c,c',l}=(c+c')a-4kl\ge 
(c+c')(a-2|k|)\ge 0\,,\end{equation} and thus
$H$ is nef.  It is well known that a nef divisor $A$ on a projective manifold
$X$ is big if and only if $A^{\dim X}>0$ (see for example, \cite[pp.\
146--147]{ShSo}).  Under the condition (\ref{nef}) we have $H^2 = 2a^2-4k^2
\ge a^2 > 0$ unless $a=k=0,$ and thus $H$ is big unless $H \equiv 0.$

\smallskip
(c): If $H$ is ample, then (\ref{a}) follows from (b) and the fact that the
ample cone is the interior of the nef cone (or simply by noting that $H\cdot
D_{1,\pm 1}>0$).  Conversely, assume (\ref{a}).  Then as we observed in (b),
$H^2>0$. Now let $\G\ne 0$ be effective, with $\G^\Sigma\equiv D_{c,c',l}$. 
Then by (\ref{a}) and (\ref{HGamma}), $H\cdot\G \ge
c+c'>0$. The ampleness of $H$ now follows by the Nakai-Moishezon
criterion.

\smallskip
(d): Assume that the conditions in (d) hold. By (c), $H$ is ample.  To show
that $H$ is globally generated, represent $H = K + L$ where $K \equiv
D_{2,\,0}$ and  $L  \equiv D_{a-2,\,k}.$ Then, in virtue of (\ref{gg})
and (b), the divisor $L$ is nef.  Furthermore,  $L^2=2(a-2)^2-4k^2\ge
(a-2)^2\ge 9$ if $a\geq 5$ and $L^2=8$ if $(a,k)=(4,0).$  So, by Reider's
Theorem, $H$ is globally generated unless there is an effective divisor
$\G=D_{c,c',l}+\G^\perp$ such that $(i)$ or $(ii)$ holds. First suppose we
have $(ii)$; by the Hodge Index Theorem, $(\G^\perp)^2 \le 0$ 
and so,
\begin{equation}\label{index} \G^2 = D_{c,c',l}^2+(\G^\perp)^2\le
D_{c,c',l}^2 =2(cc'-2l^2)\,.\end{equation} Since  $\G^2=0$ by $(ii)$,
it follows from (\ref{index}) that 
\begin{equation}\label{by-ii} cc'\ge 2l^2\,.\end{equation} Furthermore by
$(ii)$, \begin{equation}\label{lgamma} 1= L\cdot \G 
= L\cdot \G^{\Sigma}= D_{a-2,k}\cdot D_{c,c',l}
= (a-2)(c+c') -4kl\,.\end{equation} Then by (\ref{lgamma}), $kl>0$.  Hence, by
(\ref{gg}), \begin{equation}\label{kl} 4kl = 4|k||l| \le 2|l|(a-2)\,,
\end{equation} and therefore by (\ref{lgamma})--(\ref{kl}),
\begin{equation}\label{1} 1\geq (a-2)(c+c'-2|l|)\ge 2(c+c'-2|l|)\,.
\end{equation} Thus \begin{equation}\label{2l}
2|l|\ge c+c'-\half\,,\end{equation} and hence by (\ref{by-ii}) and 
(\ref{2l}), $$\left(2|l|+\half\right)^2\ge (c+c')^2\ge 4cc'\ge 8l^2\,,$$
or $2l^2-|l|\le {1\over 8}$.  Since $4l\in\Z$, it follows that
$|l|\le\half$, and then by (\ref{2l}),  $c+c'\le 1$.  According
to the claim, the only such possibility is $\G\equiv E$ or $F$, but in that
case, $L\cdot\G=a-2\ge 2$, contradicting $(ii)$.

Next suppose that $(i)$ holds.  This time, since $\G^2=-1$, (\ref{index})
yields \begin{equation}\label{by-i} cc'\ge 2l^2-\half\,.\end{equation} 
Also by $(i)$, we have 
\begin{equation}\label{lgamma'} 0= L\cdot \G = L\cdot \G^{\Sigma}= 
D_{a-2,k}\cdot D_{c,c',l}
= (a-2)(c+c') -4kl\,.\end{equation} Now, as in case $(ii)$, using
(\ref{lgamma'}) in place of (\ref{lgamma}), we obtain
\begin{equation}\label{1'} 0\geq (a-2)(c+c'-2|l|)\ge 2(c+c'-2|l|)\,.
\end{equation} and hence \begin{equation}\label{2l'}
2|l|\ge c+c'\,.\end{equation}  This time by (\ref{by-i}) and (\ref{2l'}),
$$4l^2\ge c^2+2cc'+c'{}^2\ge c^2+c'{}^2 + 4l^2-1$$ and thus $c^2+c'{}^2\le 1$.
Again the only possibility is $\G\equiv E$ or $F$, but that contradicts
$(i)$.  Therefore, by Reider's Theorem, $H$ is globally generated. The
non-speciality of $H$ follows from the Ramanujam (or Kawamata-Viehweg)
Vanishing Theorem (see for example, \cite[Ch.~VII]{ShSo}) applied to $L$.

\smallskip
(e): Suppose $H$ is a very 
ample divisor. Then the restrictions 
$H\,|\,E$ and $H\,|\,F$
are very ample, too, and hence, 
$a =  \deg H\,|\,E=\deg H\,|\,F \ge 5$. 
(Indeed, there is no very ample divisor on $C$ of degree $\le 4,$
because the genus of a smooth non-plane curve in 
$\PP^3$ of degree $\le 4$ is at most $1$.)
Furthermore, 
$2a - 4k = \deg H\,|\,D(i) \ge 5$ and 
$2a + 4k = \deg H\,|\,\Delta \ge 5,$ and so, $4|k| \le 2a-6,$
which gives (\ref{va}). 

Conversely, suppose $a \ge 5$ and 
(\ref{va}) is fulfilled. 
As above, represent $H = K + L$.
In virtue of (\ref{va}), we have
$ 2|k| \le (a-2)-1$ and so, by (c), 
$L$ is an ample divisor. We also have: 
$$L^2 = D_{a-2,\,k}^2 = 2(a-2)^2 - 4k^2
\ge 2(a-2)^2 - (a-3)^2 = a^2-2a-1 \ge 14\,,$$ 
since  $a \ge 5.$

Hence by Reider's Theorem, $H$ is very ample unless there is an effective
divisor $\G= D_{c,c',l}+\G^\perp$ such that $(i'),\ (ii')$ or $(iii')$ holds. 
However, $(i')$ cannot hold since $L$ is ample. Now suppose that $(ii')$
holds.  Since in this case, $\G^2\ge -1$, it follows from (\ref{index}) as in
case $(i)$ of (d) that (\ref{by-i}) again holds. Also by
$(ii')$, we again have (\ref{lgamma}), and repeating the argument in case
$(ii)$ of (d), we obtain $$1\ge (a-2)(c+c'-2|l|)\ge 3(c+c'-2|l|)\,.$$ Since
$4l\in\Z$, we again obtain (\ref{2l'}). However, we showed in (d)
that (\ref{by-i}) and (\ref{2l'}) imply that $\G\equiv E$ or $F$, and this
contradicts $(ii')$.  

Finally, suppose that $(iii')$ holds. Since in this case, $\G^2=0$, the
Hodge Index Theorem again yields (\ref{by-ii}) as in (d).  This time
\begin{equation}\label{lgammaiii} 2=L\cdot \G
=  L\cdot \G^{\Sigma}= (a-2)(c+c') -4kl\,.\end{equation}
Using (\ref{lgammaiii}) as in the proof of (d), we obtain 
\begin{equation}\label{2} 2\ge (a-2)(c+c'-2|l|)\ge 3
(c+c'-2|l|)\,.\end{equation}
Since $4l\in \Z$, (\ref{2}) yields (\ref{2l}).  In (d), we showed that
(\ref{by-ii}) and (\ref{2l}) imply that $\G\equiv E$ or $F$, and hence
$L\cdot\G=a-2\ge 3$, contradicting $(iii')$.  Therefore by Reider's Theorem,
$H$ is very ample. \qed

\medskip
\noindent\begin{rem} 
\label{(e)} It follows from the proof of
parts (c) and (e) of Theorem~\ref{g2ample} that if a divisor 
$H\equiv D_{a,a',k} \in \Sigma$ is ample, then 
$$2k^2 < aa' \quad \mbox{and}\quad 4|k| \le a+a'-1\,;$$
if $H$ is very ample, then $a,a'\ge 5$ and
$$ 4|k| \le a+a'-5\,. $$
\end{rem}

\begin{cor}\label{g2'ample} Let $V' = C_2$  be the symmetric square of a
genus $2$ curve $C$.  Then a divisor $H \equiv D'_{a,k} \in
\Sigma'\subset\NS(V')$ is quasi-effective, resp.\ nef, resp.\ ample, iff
the inequality  $(\ref{effective})$, resp.\ $(\ref{nef})$, resp.\
$(\ref{a})$, holds.\end{cor}

\noindent\proof This follows by Kouvidakis \cite{Kou} 
in the case where $C$ is a generic curve, and by Theorem \ref{g2ample} and
Proposition \ref{elementary} in the general case. \qed 

\medskip
Corresponding (less precise) descriptions of the globally 
generated and very
ample divisors on $V'$ are given in Theorem \ref{g2ample'} below. 
The next corollary provides some
further consequences of Theorem \ref{g2ample}.

\begin{cor} \label{min}  As above, let $V = C\times C$ 
where $C$ is a genus $2$ curve.
Then the following statements hold.

\smallskip

\noindent {\rm (a)} If $H \equiv D_{a,k} \in  \Sigma^\sym \subset \NS(V)$ 
is a very ample divisor, then the degree $n = H^2$ 
of the projective embedding $V\hookrightarrow \PP^N$ 
defined by $H$ is at least $a^2+6a-9\ge 46$.
This lower bound is achieved by the very ample divisor 
$H = 6(E+F)-\Delta\equiv D_{5,1} \in \Sigma^\sym.$

\medskip

\noindent {\rm (b)} Let 
$H \equiv D_{a,k} \in  \Sigma^\sym\subset \NS(V), \,\,\,k\neq 0,$ 
be an ample divisor. Then the divisor $2H$ is globally generated, 
and $3H$ is very ample\footnote{Recall 
\cite[II.6.1, III.7]{Mu} that if $A$ is a simple abelian variety, 
then any effective divisor $H$ on
$A$ is ample, $2H$ is globally generated and $3H$ is very ample.}. 
More precisely,
let $m_0 > 0$ be such that the divisor $m_0H$ 
is very ample but $(m_0-1)H$ is not. 
If $2|k|=a-1$ then $m_0= 3,$ and 
if $2|k|=a-2,$ then $m_0= 2.$

\medskip

\noindent {\rm (c)} Very ampleness
of a divisor $H  \equiv D_{a,k}$ 
on the surface $V$ is a numerical condition. 
Moreover, if $H= K + L\equiv D_{a,k}$, $a\ge 5$, then $H$ is ample and
globally generated if $L$ is nef, and $H$ is 
very ample if and only if $L$ is ample. 

\medskip

\noindent {\rm (d)} There is no embedding 
$V \hookrightarrow \PP^4$ 
defined by a symmetric divisor 
$H\equiv D_{a,k} \in  \Sigma^\sym\subset \NS(V)$; in particular, if $C$ has
general moduli  so that $\Sigma=\NS(V)$, then there is no
embedding $V \hookrightarrow \PP^4$  in which the images of the
generators  $E$ and $F$ are of the same degree.  
\end{cor}

\medskip

\noindent \proof (a) By (\ref{va}) we have:
$$n=H^2 = D_{a,k}^2 = 2a^2 - 4k^2 \ge 2a^2 - 
(a-3)^2 = a^2 + 6a - 9\ge 46\,,$$
since $a \ge 5$. By Theorem \ref{g2ample}, the divisor 
$H = 6E + 6F - \Delta\equiv D_{5,1}$ is very ample. 
It defines a projective embedding of $V$ of degree 
$n = 46.$

\medskip

\noindent (b) The statement follows from Theorem
\ref{g2ample}. 

\medskip

\noindent (c) Theorem \ref{g2ample}(e) tells us that very ampleness of
a divisor in $\Sigma^\sym$ is numerical.  Inequalities
(\ref{nef})--(\ref{va})
immediately imply the next statement. 
\medskip

\noindent (d) Suppose that $H\equiv D_{a,k} \in \NS(V)$ 
is a very ample divisor which defines 
an embedding $V \hookrightarrow \PP^4.$ Then $a\ge 5,$ and 
by (\ref{P4})
we get 
\begin{equation}\label{2PP4} n^2-10n=20a + 4\,.\end{equation}
{From}  (\ref{2PP4}) we obtain 
$$n=2(a^2-2k^2) \le \sqrt{20a+29}+5 <4a < a^2\,,$$
where the last two inequalities are consequences of the bound $a\ge 5.$
Hence $a^2 < 4k^2$, which contradicts (\ref{va}). \qed

\medskip

\begin{ex} A generic linear system $g_3^1$ 
on a genus 2 curve $C$ gives rise to a symmetric 
correspondence, say, $T$ on $C,$ 
which represents an effective divisor $T \equiv D_{2,1}$ on 
$V = C \times C.$ By Theorem \ref{g2ample}(b),
this divisor is nef and big, but not ample; indeed, 
$T\cdot D(i) =0.$\end{ex}

\medskip

\begin{rem} Suppose  $H\equiv D_{a,k}\in \NS(V)$, where $V$ is as in
Theorem \ref{g2ample},  $a\ge 5$ and (\ref{gg}) holds, so that by
part (d) of the theorem, $H$ is ample and globally generated.  Then we
easily see that $H$ fails to be very ample iff $2|k|=a-2$ (and $a\ge 6$) iff
equation $(i')$ of Reider's Theorem
holds with $\G = D(i)$ or $\Delta$  (the former when $k > 0$, the latter
when $k <0$). By the proof of Theorem
\ref{g2ample}(e), neither $(ii')$ nor $(iii')$ of Reider's Theorem
can hold for such $H$. Furthermore, $(i')$ can hold only for the above
choices of $\G$. (Indeed, $(i')$ implies that, in the notation from the
proof of Theorem \ref{g2ample}, $c+c' \le 2|l|$, $cc'\ge 2l^2-1$, and hence
$(c,c',l)=(1,0,\pm\half),\;(0,1,\pm\half)$ or $(1,1,\pm 1)$. By
(\ref{qeff2}), the first two cases cannot occur and the last happens only
when $\G = D(i)$ or $\Delta$.)\end{rem}

\medskip

\begin{rem} It  is well known that the correspondences on a curve $C$  form
a ring with the unit  $\Delta \equiv D_{1,-1}.$ For $g(C)=2$,
multiplication by  $D(i)\equiv  D_{1,1}$ is an involution in this
ring. It
follows that the effective cone in  the Neron-Severi group $\NS(C\times C)$
is invariant under this transform. Therefore, the same holds   for
the ample and the nef cones, as well. Since multiplication by $D(i)$ 
transforms the divisor class $D_{a,a',k}$ into $D_{a,a', -k}$, we see
why the inequalities in Theorem \ref{g2ample} are symmetric with respect
to the sign change of $k$. 
Clearly, the same is true for any hyperelliptic curve. \end{rem}

\medskip

We now consider the case where $C$ is a non-hyperelliptic smooth curve of
genus $g=3$ (as in Example~\ref{genus3}) and we let $V=C\times C$. Again we
assume that $C$ has general periods; in particular, we assume that
the Neron-Severi group $\NS(V) =\Sigma$ is generated by the classes of
$E,\,F,\,\Delta.$

\medskip

\begin{lem} \label{4.1} The class $D_{10,6} \in \NS(V)$ contains 
a unique irreducible effective curve $B.$ \end{lem}

\medskip

\noindent \proof
Consider the canonical embedding $C\hookrightarrow\PP^2$ of degree $4.$ 
Denote by $B$ the symmetric correspondence on $C$ given as the
closure in $V= C\times C$ of the set of pairs $(p,q)$ where 
$\{p,q\}=L_x\cap C\setminus\{x\}$ 
for the line $L_x\in C^*$ tangent to $C$ at some non-flex
point $x\in C$ different from $p.$ 
Write $B\equiv D_{a,k}.$  To compute $a,k$, we first note
that by the Class Formula, the dual curve $C^*$ has degree $2(g+d-1)=12.$ 
The pencil of lines through a generic point $p\in C$ represents the line 
$H=p^*$ 
in the dual projective plane $\PP^{2*}$
that is tangent to $C^*$ at the point $L_p^*\in \PP^{2*}.$  By Bezout's
theorem, $H\cap C^*$ consists of $12$ points, including the point $L_p^*$ 
of
multiplicity $2.$ Thus, there are $10$ lines through $p$ tangent to $C,$
excluding the tangent $L_p,$ and hence there are $10$ choices of
points $q$ with $(p,q)\in B$; i.e., $a=B\cdot F=10.$  To compute $k,$ we
recall that the smooth quartic $C$ has $28$ bitangent lines, and thus
$D\cdot\Delta =56,$ since each bitangent gives two points of $D\cap\Delta.$
{From} the equality $56=D\cdot\Delta=2a+2gk=20+6k$ we obtain $k=6.$  
Therefore, $$B\equiv D_{10,6}\,.$$  
In fact, $B$ is an irreducible curve. To see this, consider
the  $2:1$ map $\varphi:B\to C$ given by $\varphi(p,q)=x.$  
If $B$ were
reducible, then $B=B_1\cup B_2$ and 
$\varphi^{-1}$ would have global branches
$\psi_j:C\to B_j,$ $j=1,2.$ The projection to the
first factor $\pi_1:B\to C$ of 
degree $10$ has simple critical
points over the $24$ flexes of $C.$ On the other hand, 
by the Riemann-Hurwitz Formula, the composition 
$\pi_1\circ \psi_j\,:\,C \to C,\,\,j=1,2,$
must be an isomorphism, a contradiction. 

Since $B^2 < 0,$ $B$ is the only effective divisor in the 
numerical class $D_{10,6}.$
\qed

\medskip

\begin{rem} Some other natural correspondences that we do not use here
are the tangent correspondence 
$T\equiv D_{2,10,2}$ given by the set of pairs $(x,q)$
in the above construction, its inverse $T^{-1}\equiv D_{10,2,2}$, and the
correspondence $G\equiv D_{3,1}\equiv {1\over 4}(T+T^{-1})$ given by the
closure of the set of pairs of distinct points in the same fiber of a linear
projection $C\to\PP^1$ (see e.g. \cite{GriHarr}).

The above correspondence $B$ can be 
generalized to 
higher genera in two different ways. 
First of all, we may define it in the same 
way as above for a generic plane nodal 
curve of degree $d$ and of genus $g;$ 
then we get $B \equiv 2D_{a,k}$ where
$a=(d-3)(d+g-2)$ and $k=d+g-4.$ It is easily seen
that $B^2 < 0$ only for 
$g=3,\,\,d=4.$ 

On the other hand, following 
a suggestion by C. Ciliberto, 
for a genus $g$ curve with general moduli 
we can consider the correspondence 
$$B=\{(p,\,q) \in C \times C\,|\,p+q+(g-1)r + D \sim K_C \quad{\rm for\,\,\,some}\quad r\in C,\,\,D 
\in {\rm Div}\,(C),\,\,D \ge 0\}$$ 
(geometrically, that means that
$p$ and $q$ lie on a cut of the canonical model 
$\varphi_K(C) \subset \PP^{g-1}$ of $C$ 
by the highest osculating hyperplane say, $H_r$ at 
some other point $r \in \varphi_K(C)$).
Using de Jonqu\`eres' formula \cite[VIII.5]{ACGH} 
one can verify that  $B \equiv (g-1)(g-2)D_{a,k}$ with 
$a=g^2-g-1$ and $k=g,$ and once again, $B^2 > 0$ for any $g\ge 4$. \end{rem}

The next proposition (in the case of a generic curve)
also follows from  Kouvidakis' Theorem \ref{Kou} and Proposition
\ref{elementary}. For the convenience of the reader, we provide a
direct proof below.

\medskip

\begin{prop} \label{pr3} Let C be a 
non-hyperelliptic genus 3 curve.
Then a class $D_{a,k} \in \Sigma^\sym\subset\NS(C\times C)$ is

\noindent {\rm (a)} quasi-effective
if and only if  \begin{equation}\label{effective3} - a \le k
\le {3\over 5}a\,\,\,\,\,\,\,\,\,\,\,\,\,\,\, \Longleftrightarrow
\,\,\,\,\,\,\,\,\,\,\,\,\,\,\,\vert 
5k+a \vert \le 4a\,;\end{equation}

\noindent {\rm (b)} nef
if and only if  \begin{equation}\label{nef3} -{1\over 3} a \le k
\le {5\over 9}a\,\,\,\,\,\,\,\,\,\,\,\,\,\,\, \Longleftrightarrow
\,\,\,\,\,\,\,\,\,\,\,\,\,\,\,\vert 
9k-a \vert \le 4a\,;\end{equation}

\noindent {\rm (c)} ample if and only if
\begin{equation}\label{ample3} -{1\over 3} a < k < {5\over 9}a
\,\,\,\,\,\,\,\,\,\,\,\,\,\,\,
\Longleftrightarrow \,\,\,\,\,\,\,\,\,\,\,\,\,\,\,\vert 9k-a 
\vert < 4a\,.\end{equation} \end{prop}

\noindent \proof We first show (c): Suppose
that $D\equiv D_{a,k}$ is ample.   
Then 
\begin{equation}\label{dotdelta} D\cdot\Delta=D_{a,k}\cdot D_{1,-1} =2(a+3k) 
>0\,,\end{equation}
\begin{equation}\label{dotb} D\cdot B=D_{a,k}\cdot D_{10,6} =4(5a-9k) 
>0\,,\end{equation}
which yields  (\ref{ample3}).
Conversely, suppose that (\ref{ample3}) 
holds, and let $D\equiv D_{a,k}.$  By (\ref{ample3}),
$D^2= 2(a^2-3k^2)>0.$ Thus by the Nakai-Moishezon
criterion, it suffices to show that 
$D\cdot \Gamma >0$ for
all effective curves $\Gamma$ on $V.$ As above, 
in virtue of  (\ref{ample3}), 
we have $D\cdot \Delta > 0$ and 
$D\cdot B > 0.$
Let $\Gamma$
be an irreducible effective divisor different from  $B$ and $\Delta,$ and
let as above, $\G=\G^{\Sigma}+\G^{\perp}$ where  $\G^{\Sigma}\equiv
D_{c,c',l}\in \Sigma\otimes \Q$ and $\G^{\perp}=\G-\G^{\Sigma} \perp
\Sigma\otimes \Q$.  Then  $\G\cdot B =\G^{\Sigma}\cdot B= 10s - 36l \ge 0$
and $\G\cdot \Delta =\G^{\Sigma}\cdot\Delta= s + 6l \ge 0$  where $s:=c+c',$
i.e.  \begin{equation}\label{ls} -{1\over 6} s \le l \le {5 \over 18}
s.\end{equation}  {From}  (\ref{ample3})
and (\ref{ls}), we obtain  
$kl \le |k||l|<  {25\over 162} as.$
Furthermore, since by  (\ref{ample3}), 
$a>0$, and clearly $s=c+c'>0,$ 
we have
$$D\cdot\G= D\cdot\G^{\Sigma}=as-6kl >{2\over 27}as >0\,.$$ 
Thus, $D$ is ample.

To show (b), we proceed exactly as above, first noting that if $D$ is
nef, then the inequalities $D\cdot\D\ge 0,\ D\cdot B\ge 0$ yield
(\ref{nef3}).  Conversely, (\ref{nef3}) implies that $D\cdot\D\ge 0,\
D\cdot B\ge 0$, and $D\cdot\G>0$ for irreducible curves $\G$ different
from $B$ and $\D$.

We now show (a):  First suppose that $k\ge 0.$  If
$D_{a,k}$ is quasi-effective, then by (\ref{nef3}), we have $D_{a,k}\cdot
D_{c,l} = 2(ac-3kl) \ge 0$ whenever
$l={5\over 9}c>0$, and thus $k\le {3\over 5}a.$  Conversely, if $0\le
k\le {3\over 5}a$, then $D_{a,k}\equiv {1\over3}D_{3a- 5k ,0} + {k\over
6}B$ is quasi-effective. Now suppose that $k\le 0.$ If $D_{a,k}$ is
quasi-effective then $ac-3kl \ge 0$ for $l=-{1\over 3}c<0$, 
which gives $k\ge -a.$  Conversely,
if $-a\le k\le 0$, then $D_{a,k}\equiv D_{a+k,0} +|k|\D$ is
quasi-effective. \qed

\medskip

To simplify the discussion, we assume 
in the sequel that the curve $C$ has
general moduli so that $\NS(V)=\Sigma$, resp.\ $\NS(V')=\Sigma'$.

\medskip

\begin{thm}\label{veryample3} Let C be a 
non-hyperelliptic genus $3$ curve with 
$\NS(V) = \Sigma$. Then a divisor 
$H\equiv D_{a,k} \in \Sigma^\sym\subset \NS(V)$ is very ample
if    $a \ge 7$,
$(a,k)\ne (8,2)$ and
\begin{equation}\label{va3} -{1\over 3} (a-4) \le k
\le {5\over 9}(a-4)\,;\end{equation} 
$H$ is ample and globally generated for 
the additional values $(6,0)$ 
and $(8,2)$ of $(a,k)$. In all these cases $H$ is non-special. \end{thm}

\noindent \proof Suppose $a\ge 8$ and (\ref{va3}) holds;
represent $H = K + L$, where $K=K_V\equiv D_{4,0}$ and $L=H - K\equiv
D_{a-4,k}.$ By Proposition~\ref{pr3}, the divisor $L$ is nef.
If $a\ge 13$, then by the inequalities (\ref{va3}), we have
$$L^2 = 2\left((a-4)^2-3k^2\right) \ge {4\over 27}(a-4)^2\ge 12\,.$$ 
One can easily
check that also $L^2\ge 12$ 
for all the integer solutions of (\ref{va3}) with 
$7\le a\le 12$, except for $(a,k)=(8,2)$.
Thus, Reider's criterion can be
applied; that is, the divisor $H=K+L$ is very ample unless there exists an
effective divisor $\G\equiv D_{c,c',l}$ such that
\begin{equation}\label{ii'3}
\G^2=2(cc'-3l^2)=0 \quad {\rm and}\quad 1 
\le L\cdot \G = (a-4)s - 6kl \le 2\end{equation} 
or 
\begin{equation}\label{i'3} \G^2=2(cc'-3l^2)=-2 \quad {\rm
and}\quad L\cdot \G = (a-4)s - 6kl =0\,,\end{equation} 
where $c,c'\ge 0$, $s=c+c'>0.$ (The
fact that $\G^2$ is even eliminates 
the other possibilities in Reider's
Theorem.) 

We first consider the case $a=7$.  By (\ref{va3}), $|k|\le 1$ when
$a=7$. If $k=0$, then $L\cdot \G=3s$ and thus neither (\ref{ii'3}) nor
(\ref{i'3}) can hold. Now suppose that $|k|=1$, so that $L\cdot \G=3s\pm 6
l$.  Hence (\ref{ii'3}) cannot hold; suppose further that (\ref{i'3}) is
satisfied, i.e., $$ cc'=3 l^2-1 \qquad {\rm and} \qquad s=2|l|\,.$$  But
then,
$$4l^2=s^2=(c+c')^2\ge4cc'=12l^2-4$$ or $2l^2\le 1$. This implies that
$l=s=0$, a contradiction, so (\ref{i'3}) also cannot hold. Thus in the
sequel we assume that $a \ge 8.$

First suppose there is an effective divisor $\G$
satisfying (\ref{ii'3}).  The inequality $(a-4)s - 6kl \le 2$
implies that $kl>0$, and hence from (\ref{va3}) and (\ref{ii'3})
we get $${10\over 3}(a-4)|l| \ge
6|k||l|= 6kl \ge (a-4)s-2$$ or
$$3s-10|l| \le {6\over a-4}\;.$$
Since $a\ge 8$, this yields
\begin{equation} \label{sle} 3s-10|l|\le 1\;. \end{equation} 
By (\ref{ii'3}),
\begin{equation} \label{sge} s^2= (c+c')^2\ge 4cc'=12 l^2\,.\end{equation}
Combining (\ref{sle}) and (\ref{sge}), we get
$${25\over 3}s^2\ge 100 l^2 \ge  (3s-1)^2\,,$$
or $$2s^2-18s+3\le 0\,.$$

Hence $s\le 8$ and thus by (\ref{sge}) and the fact that $kl>0$, we
have $1\le |l|\le 2$. If $l=\pm 1$, then $cc'=3$ and thus
$c=1,\ c'=3$ (or vice versa)  so that $s=4$.  But then $3s-10|l|=2$,
contradicting (\ref{sle}).  

Now suppose that $l=-2$.  Then $k<0$ and hence by (\ref{va3}) we have
$3|k|\le a-4$. Thus by (\ref{ii'3}), $$3|k|(s-4) = 3|k|s
-12|k| \le (a-4)s-12|k|\le 2\,.$$ Therefore $s\le
4$, which contradicts the fact that by (\ref{sge}), $s^2\ge 48$.

It remains to consider $l=2$.  In this case, (\ref{sle}) yields $s\le
7$. Let $$\G=\sum_{j=1}^m \G_j\equiv \sum_{j=1}^m D_{c_j,c'_j,l_j}$$ be the
decomposition of $\G$ into irreducible effective divisors. Note that
$c_j,c'_j \ge 0$, $\sum  (c_j +c'_j)=s\le 7$, $\sum l_j = l=2$, and thus the
$\G_j$ are all distinct from $B\equiv D_{10,10,6}.$  
Therefore,  as in the proof
of Proposition~\ref{pr3} (see (\ref{ls})), $B\cdot \G_j \ge 0$ gives 
\begin{equation}\label {lj} l_j\le {5\over 18}
(c_j+c'_j)\,.\end{equation}  Summing (\ref{lj}) over $j$ yields $2\le
{5\over 18}\cdot 7$, a contradiction. 

Next, suppose there is an effective divisor $\G$
satisfying (\ref{i'3}).  As before, we conclude that $kl>0$. Hence from
(\ref{va3}) and (\ref{i'3}) we get $${10\over 3}(a-4)|l| \ge
6|k||l|= 6kl= (a-4)s$$ or
\begin{equation} \label{sle'} s\le {10\over 3}|l|\;.\end{equation} 
By (\ref{i'3}), \begin{equation}\label{3l}
3l^2=cc'+1 \end{equation} and hence by (\ref{sle'}),
\begin{equation} \label{ne} {100\over
9}l^2\ge  s^2=(c+c')^2\ge 4cc'=12 l^2 -4\,.\end{equation} By (\ref{ne}) and the
fact that $kl>0$, we have again $1\le |l|\le 2$.  If $|l|=2$, then (\ref{3l})
yields $cc'=11$, which contradicts the fact that by (\ref{sle'}), $c+c'=s\le
6$. 

If $|l|=1$, we have $cc'=2$ and $s=3$; 
therefore $\G\equiv D_{2,1,\pm
1}$ or $D_{1,2,\pm 1}.$ If $\G\equiv D_{2,1,-1}$ or $D_{1,2,-1}$, then
(\ref{i'3}) yields $k=-\half (a-4)$, which contradicts (\ref{va3}). On the
other hand, if $\G\equiv D_{2,1,1}$ or $D_{1,2,1}$, then we again
decompose $\G$ into irreducible effective divisors: $$\G=\sum_{j=1}^m
\G_j\equiv \sum_{j=1}^m D_{c_j,c'_j,l_j}$$ Since $\sum  (c_j +c'_j)=3$, $\sum
l_j = 1$, the $\G_j$ must all be distinct from $B\equiv D_{10,10,6}$ and
thus (\ref{lj}) holds.  Summing (\ref{lj}) over $j$ then yields
$1 \le {5\over 18}\cdot 3$, a contradiction.  To summarize,
$\G$ cannot satisfy (\ref{ii'3}) or (\ref{i'3}); therefore $H$ is very ample
by Reider's Theorem \ref{re}.

To obtain the second statement, suppose $H\equiv D_{6,0}$, resp.\
$D_{8,2},$ and write $H=K+L$ as before. Then in both cases, $L^2=8$, and
(\ref{va3}) holds, so by Proposition
\ref{pr3}(b), $L$ is nef.  For any divisor
$\G\equiv D_{c,c',l}\in\NS(V)$, we have $\G^2\neq  -1$ and $L\cdot\G =
2s$, resp.\ $4s-12 l$, and hence $L\cdot\G \neq 1$.
Thus, $H$ is globally generated by Reider's Theorem; furthermore,
$H$ is ample by Proposition \ref{pr3}(c). As above, the non-speciality of $H$ 
follows from the Ramanujam Vanishing Theorem. \qed 

\medskip

\begin{rem} \label{rm4.3} 
The smallest possible degree 
of a projective embedding $V \hookrightarrow \PP^N$ provided
by  the above
theorem is $92,$ given by the very ample divisor
$D_{7,1}$. Recall (Example~\ref{genus3}) that the canonical
divisor $K_C\equiv D_{4,0}$  is also very ample and
gives an embedding of degree 32. \end{rem}

\medskip

Applying the same methods as in Theorems \ref{g2ample}(d)--(e) and
\ref{veryample3} above, we describe in the next theorem some of the
non-special globally generated, resp.\ very ample, divisors on the symmetric
square $V'=C_2$ of generic curves $C$ of genera 2, 3 and 4. (Recall that if $C$
has genus 2, then $C_2$ is an abelian surface with a point blown up and is
consequently non-hyperbolic.)

\medskip

\begin{thm} \label{g2ample'}  Let $V' = C_2$ 
where $C$ is a genus $g\ge 2$ curve 
with general moduli, so that $\NS(V') = \Sigma'.$ 
Let  $H \equiv D'_{a,k} \in \NS(V')$ be a divisor.

\medskip

\noindent For $g=2$ the divisor $H$ is

\smallskip\noindent {\rm (a)} non-special, nef and globally generated if  
$a\ge 4$ and
\begin{equation}\label{g2g'} 2|k-1| \le a-2\,;\end{equation}

\smallskip\noindent {\rm (b)} very ample if $a \ge 5$ and 
\begin{equation}\label{v2a'} 2|k-1| \le a-3\,.\end{equation}
On the other hand, if $H$ is very ample, then $a\ge 5$ and
\begin{equation}\label{v2a''} 1-a \le 2(k-1) \le a-3\,.\end{equation}

\medskip

\noindent For $g=3$ the divisor $H$ is

\smallskip\noindent {\rm (a$'$)} non-special, ample and globally generated if  
$a\ge 6$, $(a,k)\ne(7,3)$ and
\begin{equation}\label{gg3} 3-a \le 3(k-1) \le {5\over
3}a-5\,;\end{equation} 

\smallskip\noindent {\rm (b$'$)}
very ample if $a \ge 7$, $(a,k)\ne (7,3)$ or $(9,4)$,
and  
\begin{equation}\label{va3'} 4-a \le 3(k-1) \le {5a-16\over 3}
\,.\end{equation}  
On the other hand, if $H$ is very ample, then $a\ge
4$ and 
\begin{equation}\label{va3''} -1-a \le 3(k-1) \le
{5a-11\over 3}\,.\end{equation}

\medskip

\noindent For $g=4$ the divisor $H$ is

\smallskip\noindent {\rm (a$''$)} non-special, ample and globally generated
if  $a \ge 8$ and
\begin{equation}\label{gg4} 9-a \le 4k \le 2a-8\,;
\end{equation}

\smallskip\noindent {\rm (b$''$)}
very ample if $a \ge 9$ and
\begin{equation}\label{va4} 9-a \le 4k \le  2a-10\,.
\end{equation}
On the other hand, if $H$ is very ample, then $a\ge 6$ and
\begin{equation}\label{va4''} 3-a \le 4k\le 2a-6\,.\end{equation}
\end{thm}

\noindent \proof (a): Let $g=g(C)=2.$ 
Set $L=H-K\equiv D'_{a-1,k-1}.$ By Theorem
\ref{Kou}, in view of (\ref{g2g'}) the divisor $L$ is ample. 
Inequality (\ref{g2g'}) implies that $2|k| \le a,$ and hence by Theorem
\ref{Kou}, $H$ is nef. 
By the Kodaira Vanishing Theorem,
$H$ is also non-special. 

Since $a\ge 4$, we have
$$L^2 = (a-1)^2-2(k-1)^2 \ge {a^2\over 2} -1 \ge  7\,.$$ 
Thus by Reider's Theorem, $H$ is globally generated unless 
one of the cases $(i)$ or $(ii)$ of this theorem happens. 
But $(i)$ is impossible since $L$ is ample, and $(ii)$ is impossible since
there is no non-zero divisor $\G\equiv D'_{c,l}$ on $V'$ with 
$\G^2 = c^2-2l^2 =0.$ 
\qed

(b): Assume first that $a\ge 5$ and (\ref{v2a'}) holds. Then as above,
$L=H-K$ is an ample divisor, and $L^2\ge a^2/2 -1 \ge  11$.
Thus by Reider's Theorem, $H$ is very ample unless 
one of the cases $(i')$-$(iii')$ of this theorem happens. 
The cases $(ii')$ with $\G^2=0$, $(i')$, and $(iii')$
are excluded by the same reasons as above. Then we are left with the
possibility that 
\begin{equation}\label{biii'} \G^2=c^2-2l^2=-1\quad {\rm and}\quad L\cdot \G
= (a-1)c-2(k-1)l=1\end{equation} for an effective divisor $\G\equiv D'_{c,l}
\in \NS(V')$. By (\ref{biii'}), $(k-1)l>0$, and furthermore by
(\ref{v2a'}) and (\ref{biii'}), $$1= L\cdot\G \ge (a-1)(c-|l|)+4|l|\,;$$
hence, $|l|>c$. (If $|l|=c$, the above inequality would yield $1\ge
4|l|=4c\ge 4$, a contradiction.)  
But then $c^2-2l^2 <-c^2 \le -1$ contradicting
(\ref{biii'}). Therefore by Reider's Theorem, $H$ is very ample.

Suppose now that $H$ is very ample. Since  $E'$, resp.\
$\Delta'$, is a smooth genus 2 curve in $V',$
and the restriction $H\,|\,E'$, resp.\ $H\,|\,\Delta'$, 
is very ample, we have
that $$\deg(H\,|\,E')=H\cdot E' = a \ge 5\,,$$ 
resp. 
$$\deg(H\,|\,\Delta')=H\cdot \Delta' = 
2D'_{a,k}\cdot D'_{1,-1}=2(a+2k) \ge 5\,.$$ Therefore, $2(k-1) \ge 1-a$.
On the other hand, since $H$ is ample, by Theorem \ref{Kou}, we have $2k < a,$
or $2(k-1) \le a-3.$ Finally, from these we get the inequality (\ref{v2a''}).
\qed

\smallskip

(a$'$), (b$'$): For $g=3$ we have $K=K_{V'}\equiv D_{3,1}$ and 
$L=H-K\equiv D'_{a-3,k-1}.$ In virtue of Theorem \ref{Kou}, 
$L$ is nef, resp.\ ample, iff (\ref{gg3}), resp.\ (\ref{va3'}),
holds. The inequality (\ref{gg3}) implies that
\begin{equation}\label{3g'} -a < 3k < 5a/3\,\end{equation}
and hence (again by Theorem \ref{Kou}),  $H$ is ample. 
By Ramanujam's Vanishing Theorem, the divisor $H=K+L$
is non-special as soon as $L$ is nef. 

Furthermore,  assuming (\ref{gg3}), resp.\ (\ref{va3'}), we get the
inequality 
\begin{equation}\label{*} L^2 = (a-3)^2-3(k-1)^2 \ge {2\over 27}
(a-3)^2\,,\qquad {\rm resp.}\qquad L^2>{2\over 27} (a-3)^2\,,\end{equation}
and hence $L^2\ge 10$ if $a \ge 15$.  By checking all possible
integer values of $a,\,k$ with $a\le 14$ satisfying the conditions of
(a$'$), resp.\ (b$'$), one easily verifies that $L^2\ge 5$, resp.\ $\ge
10$,  for these values. Therefore by Reider's Theorem, under the assumptions
on $(a,\,k)$ of   (a$'$), resp.\ (b$'$), $H$ is globally generated,
resp.\ very ample, unless there is an effective divisor 
$\G\equiv D_{c,l}\in\NS(V')$ 
for which one of the conditions $(i)$--$(ii)$, resp.\ $(i')$--$(iii')$,  of
this theorem holds. Note that the diophantine equations $\G^2=c^2-3l^2=0$ and
$\G^2=c^2-3l^2=-1$ have no solutions (since the latter cannot hold modulo
3). Finally, $L\cdot \G=0$ is impossible assuming (\ref{va3'}), 
because in that case the divisor $L$ is ample.  Thus,  
$(i)$--$(ii)$, resp.\ $(i')$--$(iii')$, cannot hold.

Next we assume that the divisor $H$ is very ample. Then so are also
the restrictions 
$H\,|\,\Delta'$ and $H\,|\,B'$ where 
$B'=\pi(B)\equiv D'_{10,6}\in \NS(V').$ 
Here $\Delta'$ is a smooth reduced curve on $V'$ isomorphic to $C,$ 
and hence $\deg (H\,|\,\Delta') \ge\delta (\Delta') = \delta (C)=4,$
i.e. $$H\cdot \Delta' = 2D'_{a,k}\cdot D'_{1,-1} = 2(a+3k) \ge 4\,.$$
By the construction, the curve $B'$ is birationally equivalent 
to the curve $C;$ in particular, the geometric genus of $B'$ equals 3.
Hence, $\deg (H\,|\,B') \ge \delta(B') \ge 4,$ or
$$H\cdot B' = D'_{a,k}\cdot D'_{10,6} = 10a - 18k \ge 4.$$ 
These inequalities provide (\ref{va3''}). \qed

\smallskip

(a$''$), (b$''$): For $g=4$ we have $K = K_{V'} \equiv D_{5,1}$ and so, 
$L = H-K\equiv D_{a-5,k-1}.$ Thus by the Kouvidakis Theorem \ref{Kou}, 
the divisor $L$ is nef iff (\ref{gg4}) holds. 
Under these inequalities we have $-a < 4k < 2a,$ which implies
(again due to Theorem \ref{Kou}) that the divisor $H$ is ample. 
By the Ramanujam Vanishing Theorem, it is non-special. 

As in the genus 3 case above, we easily verify that $L^2 \ge 5$,
resp.\ $L^2 \ge 10$, under the conditions of (a$''$), resp.\ (b$''$).  Hence,
by Reider's Theorem, $H$ is globally generated resp.\ very ample,
unless (for an effective divisor $\G\equiv D_{c,l}\in \NS(V')$) one of the
conditions $(i)$--$(ii)$, resp.\ $(i')$--$(iii')$, holds. We note that
the diophantine equations $\G^2 = c^2-4l^2=-1$ and $\G^2 = c^2-4l^2=-2$ have
no solutions (since neither can hold modulo 4). Hence, $\G^2 = 0$, that is, $\G
\equiv \alpha D_{2,\pm 1}$ where $\alpha > 0$. But then we would have
$$L\cdot \G= \alpha(2(a-5)\mp 4(k-1))=1 \quad {\rm or}\quad 2\,.$$ 
Therefore,
$\alpha = 1,$ that is, $\G \equiv D_{2,\pm 1},$ and $L\cdot \G = 2$,
i.e. $(iii')$ holds. Since the cases $(i)$ and $(ii)$ have been 
eliminated, $H$ is
globally generated provided that the assumptions in (a$''$) are fulfilled.

Assume that $(iii')$ holds. Then we have
$a+2k=8$ if $\G \equiv D_{2, -1},$ and $a=2k+4$ if $\G \equiv D_{2,1}.$
In the first case (\ref{gg4}) yields $6\le a \le 7$ 
which is excluded by the assumptions of (b$''$); 
the second case contradicts (\ref{va4}). Thus $H$ is very
ample. 

To show the last statement of (b$''$), we now suppose
that $H$ is very ample. Then so are the restrictions 
$H\,|\,E',\ H\,|\,\Delta'$ and $H\,|\,B''$ where $B''\equiv
D_{2,1}\in\NS(V')$ is a smooth curve on $V'$ isomorphic to $C,$ provided by
any of the two linear pencils $g^1_3$ on the curve $C$. (Recall
\cite[IV.5.5.2]{Hart} that a generic genus 4 curve $C$ possesses  exactly
two such pencils.) Since $E' \simeq \Delta' \simeq B'' \simeq C$  and
$\delta(C)=6$ (see Remark \ref{rm2.3} above),
we have $$a=H\cdot E'\ge \delta(C) = 6\,,\ 2(a+4k) = H\cdot \Delta'
\ge \delta(C) = 6\,,\ 2a-4k=H\cdot B'' \ge \delta(C) = 6\,.$$
This proves the inequalities $a\ge 6$ and (\ref{va4''}). \qed

\end{document}